\documentclass[11pt]{article}
\usepackage{epsfig}
\begin{document}

\author{Casim Abbas \\Department of Mathematics\\Michigan State University\\Wells Hall\\East Lansing, MI 48824\\USA}  
\title{Pseudoholomorphic strips in Symplectisations III: Embedding properties and Compactness} 
\date{\today}
\maketitle


\footnote{This material is based upon work supported by the National Science Foundation under Grant No. 0196122 and by a New York University Research Challenge Fund Grant. Part of this paper was written while the author visited ETH Z\"{u}rich. He would like to thank FIM and Prof. Eduard Zehnder for their hospitality.}  
\tableofcontents

%
%
\parindent0ex
\parskip1ex plus0.4ex minus0.2ex

\newtheorem{lemma}{Lemma}[section]
\newtheorem{theorem}[lemma]{Theorem}
\newtheorem{corollary}[lemma]{Corollary}
\newtheorem{proposition}[lemma]{Proposition}
\newtheorem{definition}[lemma]{Definition}
\newtheorem{conjecture}{Conjecture}
\newtheorem{example}[lemma]{Example}
\newtheorem{condition}[lemma]{Condition}


{\catcode`"=12 \gdef\hex{"}}
\mathchardef\laplace=\hex0001
\let\rmDelta=\laplace
\mathchardef\nabla=\hex0272
\makeatletter
\def\@@dalembert#1#2{\setbox0\hbox{$#1\mathrm I$}
  \vrule height\ht0 depth\z@ width.04\ht0
  \rlap{\vrule height\ht0 depth-.96\ht0 width.8\ht0}
  \vrule height.1\ht0 depth\z@ width.8\ht0
  \vrule height\ht0 depth\z@ width.1\ht0 }
\def\dalembert{\mathbin{\mathpalette\@@dalembert{}}\,}
\makeatother
\newcommand{\qed}{\begin{flushright}$\dalembert$\end{flushright}}
\newcommand{\degree}{\mbox{deg}}
\newcommand{\ind}{\mbox{ind}}
\newcommand{\abs}{|}
\newcommand{\ve}{\varepsilon}
\newcommand{\Id}{\mbox{Id}}
\newcommand{\GL}{\mbox{GL}}
\newcommand{\eps}{\varepsilon}
\newcommand{\To}{\longrightarrow}
\newcommand{\Real}{{\bf{R}}}
\newcommand{\Complex}{{\bf{C}}}
\newcommand{\RM}{{\bf{R}}\times M}
\newcommand{\RL}{{\bf{R}}\times {\mathcal L}}
\newcommand{\tu}{\tilde{u}}
\newcommand{\pil}{\pi_{\lambda}}
\newcommand{\pas}{\partial_s}
\newcommand{\pat}{\partial_t}
\newcommand{\od}{\{0\}\times{\mathcal D}^{\ast}}
\newcommand{\SC}{{\mathcal S}}
\newcommand{\ho}{\mbox{Hom}}
\newcommand{\oo}{{\mathcal O}}

\section{Introduction, Notation, Results}
This paper continues the investigation initiated in the papers \cite{part1} and \cite{part2} aimed at developing suitable pseudoholomorphic curve techniques for the investigation of the Chord problem in three dimensional closed contact manifolds. A contact form $\lambda$ on a $(2n+1)$--dimensional manifold $M$ is a 1--form so that $\lambda\wedge(d\lambda)^n$ is never zero. The hyperplane field $\xi=\ker\lambda$ is called the contact structure associated to $\lambda$. A Legendrian submanifold ${\mathcal L}$ is an n--dimensional submanifold of $M$ such that $\lambda|_{T{\mathcal L}}\equiv 0$. There is a distinguished vector field associated with $\lambda$, called the Reeb vector field $X_{\lambda}$. It is defined by the equations 
\[
i_{X_{\lambda}}d\lambda\equiv 0\ \mbox{and}\ i_{X_{\lambda}}\lambda\equiv 1.
\]
A characteristic chord for a given Legendrian submanifold ${\mathcal L}$ is a trajectory $x(t)$ of the Reeb vector field such that $x(0),x(T)\in {\mathcal L}$ for some $T>0$. We usually also ask for $x(0)\neq x(T)$. The question is the following: Given a manifold $M$ with contact form $\lambda$ and a Legendrian submanifold ${\mathcal L}$, is there a characteristic chord ? For very special $M$ and ${\mathcal L}$ characteristic chords are known in Hamiltonian mechanics as 'brake--orbits', and they were investigated since the 1940's by H. Seifert and many others (see for example \cite{Ambrosetti}, \cite{Seifert}, \cite{Weinstein}). In the case where $M$ is the three dimensional sphere and $\xi$ is the standard tight contact structure existence of characteristic chords for any Legendrian knot was conjectured by V.I. Arnold \cite{Arnold} and proved by K. Mohnke \cite{Mohnke} after a partial result by the author \cite{A3}. The aim is to establish a method based on filling by pseudoholomorphic curves which is able to detect characteristic chords in general closed three dimensional contact manifolds $M$. The purpose of this paper together with the previous papers \cite{part1} and \cite{part2} is to establish the filling method while we apply it in the forthcoming paper \cite{part4} to obtain an existence result for characteristic chords. Let us describe the boundary value problem which we are going to investigate. We consider a three dimensional closed manifold $M$ with contact form $\lambda$. Moreover, we assume that ${\mathcal L}\subset M$ is a homologically trivial Legendrian knot, and ${\mathcal D}$ is an embedded surface bounding ${\mathcal L}$. A point $p\in {\mathcal D}$ is called singular if the contact plane $\xi_p$ at the point $p$ is identical with the tangent plane to the surface ${\mathcal D}$. We denote the set of all singular points on ${\mathcal D}$ by $\Gamma$. Generically, the set $\Gamma$ is finite. We will call a complex structure $J:\xi\rightarrow\xi$ admissible if $d\lambda\circ(\Id\times J)$ is a metric on $\xi$. We extend such an admissible $J$ to an almost complex structure $\tilde{J}$ on the symplectisation ${\bf R}\times M$ of $M$ by demanding
\[
\tilde{J}(a,u)\frac{\partial}{\partial t}=X_{\lambda}(u)\ \ \mbox{and}\ \ \tilde{J}(a,u)X_{\lambda}(u)=-\frac{\partial}{\partial \tau}
\]
for $(a,u)\in {\bf R}\times M$, where $\tau$ denotes the coordinate in the ${\bf R}$--direction. The Seifert surface ${\mathcal D}$ can be perturbed near its boundary, leaving the boundary fixed, in order to achieve a certain normal form in local coordinates (see proposition \ref{3.1.1.} below). In the earlier paper \cite{part2} we also chose a particular complex structure $J:\xi\rightarrow\xi$ near $\{0\}\times{\mathcal L}$ (see (\ref{local-formula-for-J}) below). We will also adhere to these choices in this paper. We note however, that the intersection result below (theorem \ref{intersection-of-families-0}) does not use them. In this paper we will study pseudoholomorphic strips in the symplectisation ${\bf R}\times M$ with satisfy the following mixed boundary condition:

\begin{equation}\label{main-boundary-value-problem}
\left\{\begin{array}{ll}
          \tilde{u}=(a,u):S\longrightarrow\RM & \\
          \partial_s\tu+\tilde{J}(\tu)\partial_t\tu=0 & \\
          \tu(s,0)\subset{\bf R}\times{\cal L} & \\
          \tu(s,1)\subset\{0\}\times{\cal D}^{\ast} & \\
           E(\tu)<+\infty
          \end{array}\right.
\end{equation}
where $S:={\bf R}\times[0,1]$, ${\cal D}^{\ast}$ is the
spanning surface ${\cal D}$ without the set of singular points $\Gamma$ on ${\mathcal L}=\partial{\mathcal D}$ and where $E(\tu)$ is the energy of $\tu$ defined by
\[
E(\tu):=\sup_{\phi\in\Sigma}\int_S\tu^{\ast}d(\phi\lambda)\ ,\ \Sigma:=\{\phi\in C^{\infty}({\bf R},[0,1])\,|\,\phi'\ge 0\}.
\]
Solutions to (\ref{main-boundary-value-problem}) exist locally near elliptic singular points on the boundary if the Seifert surface ${\mathcal D}$ has been perturbed into normal form as in proposition \ref{3.1.1.} and if the complex structure $J:\xi\rightarrow\xi$ has been chosen appropriately:

\begin{theorem}\label{local-existence-theorem}
Let $(M,\lambda)$ be a three dimensional contact manifold. Moreover, let ${\mathcal L}$ be a Legendrian knot which bounds an embedded surface ${\mathcal D}'$ so that the characteristic foliation has only finitely many singular points. Then there is another embedded surface ${\mathcal D}$ which is a smooth $C^0$--small perturbation of ${\mathcal D}'$ having the same boundary and the same singular points as ${\mathcal D}'$ and a $d\lambda$--compatible complex structure $J:\ker\lambda\rightarrow\ker\lambda$ so that the following is true: Near each elliptic singular point $e\in\partial{\mathcal D}={\mathcal L}$ there are embedded solutions $\tu_{\tau}$ , $0<\tau<1$ to the boundary value problem (\ref{main-boundary-value-problem}) with the properties:
\begin{itemize}
\item $\tu_{\tau}(S)\cap\tu_{\tau'}(S)=\emptyset$ if $\tau\neq\tau'$,
\item $\tu_{\tau}\rightarrow e$ uniformly with all derivatives as $\tau\rightarrow 0$,
\item the family $\tu_{\tau}$ depends smoothly on the parameter $\tau$,
\item each map $u_{\tau}$ is transverse to the Reeb vector field, i.e. $\pil\pas u_{\tau}(z)\neq 0$ for all $z\in S$,
\item The Maslov indices $\mu(\tu_{\tau})$ all equal zero.
\end{itemize}
\end{theorem}

See section \ref{local-existence-of-sol} of this paper for the existence statement and \cite{part2} for the statement about the Maslov indices. The solution family $\tu_{\tau}$ above is unique up to parametrization. One of the main results of this paper is the following compactness result: 
\begin{theorem}\label{compactness-result}
Let $(\tu_{\tau})_{0\le\tau<\tau_0}=(a_{\tau},u_{\tau})_{0\le\tau<\tau_0}$ be a smooth family of embedded solutions to the boundary value problem
\[
\begin{array}{ll}
          \tilde{u}=(a,u):S\longrightarrow\RM & \\
          \partial_s\tu+\tilde{J}(\tu)\partial_t\tu=0 & \\
          \tu(s,0)\subset{\bf R}\times{\cal L} & \\
          \tu(s,1)\subset\{0\}\times{\cal D}\backslash\Gamma & \\
          u(0,0)=e & \\
           0< E(\tu)<+\infty & 
          \end{array},
\]
where ${\mathcal D}\subset M$ is an embedded surface bounding the Legendrian knot ${\mathcal L}$, $\Gamma\subset {\mathcal D}$ is the set of singular points and $e\in\Gamma\cap{\mathcal L}$ is an elliptic singular point on the boundary of ${\mathcal D}$. We impose the following conditions on the solutions $\tu_{\tau}$:
\begin{itemize}
\item $\tu_{\tau'}(S)\cap\tu_{\tau''}(S)=\emptyset$ if $\tau'\neq \tau''$,
\item 
\[
\mbox{dist}\left(\bigcup_{0<\delta\le\tau<\tau_0}\{u_{\tau}(s,1)\,|\,s\in{\bf R}\}\,,\,\Gamma\right)>0,
\]
\item For small $\tau$ the solutions $\tu_{\tau}$ coincide with the local solutions of theorem \ref{local-existence-theorem} near $e$, 
\item There is a uniform gradient bound, i.e.
\[
\sup_{0\le\tau<\tau_0}\|\nabla\tu_{\tau}\|_{C^0(S)}<\infty.
\]
\end{itemize}
Then for every sequence $\tau'_k\nearrow\tau_0$ there is a subsequence $\tau_k$ such that the family $\tu_{\tau_k}$ converges in $C^{\infty}_{loc}$ to another solution (as $k\rightarrow\infty$) $\tu_{\tau_0}$ with finite energy such that also $\mbox{dist}\big(\{u_{\tau_0}(s,1)\,|\,s\in{\bf R}\}\,,\,\Gamma\big)>0$. Moreover
\begin{enumerate}
\item every sequence $\tau_k$ yields the same limit, i.e. $\tu_{\tau_0}=\lim_{\tau\nearrow\tau_0}\tu_{\tau}$, and the convergence is uniform on $S$ with all deerivatives,
\item $\tu_{\tau_0}$ is an embedding,
\item the Maslov--index $\mu(\tu_{\tau_0})$ of $\tu_{\tau_0}$ equals $0$,
\item The solution $\tu_{\tau_0}(s,t)$ has the same rate of decay for large $|s|$ as the maps $\tu_{\tau}$, i.e. $|\lambda_{\pm}|=\frac{\pi}{2}$ if the asymptotic formula, theorem \ref{asymptotic-formula-theorem}, is applied to $\tu_{\tau_0}$,
\item $\tu_{\tau_0}(S)\cap\tu_{\tau}(S)=\emptyset$ for all $0\le\tau<\tau_0$.
\end{enumerate}
\end{theorem}
{\bf Remarks:}\\
We have shown in the paper \cite{part1} that finiteness of energy implies the existence of the limits $\lim_{s\rightarrow\pm\infty}\tu_{\tau}(s,t)\in\{0\}\times{\mathcal L}$. See definition \ref{admissible-normal-vector} below and \cite{part2} for the definition and properties of the Maslov--index. It is a well--known fact that the uniform gradient bound implies that every sequence $(\tu_{\tau_k})_{k\in{\bf N}}$, $\tau_k\rightarrow\tau_0$  from the solution family $\tu_{\tau}$ has a $C^{\infty}_{loc}$ convergent subsequence. The theorem asserts that the whole family $\tu_{\tau}$ converges as $\tau\rightarrow\tau_0$ uniformly on $S$ and not just on compact sets. It is also not hard to see that the energy of $\tu_{\tau_0}$ is finite which implies exponential decay of the solution at infinity (theorem \ref{asymptotic-formula-theorem}). The main points of theorem \ref{compactness-result} are the following:
\begin{itemize}
\item The limit $\tu_{\tau_0}$ has the same exponential decay rate and the same Maslov--index as the rest of the family $\tu_{\tau}$ for $\tau<\tau_0$. It decays at the slowest possible rate $\sim e^{-\frac{\pi}{2}|s|}$. Note that the convergence is initially only in $C^{\infty}$ {\it on compact sets} which does not permit us to draw any conclusions about the behavior of $\tu_{\tau_0}$ at infinity from the asymptotic behavior of the solutions $\tu_{\tau}$. Additional arguments are required here.
\item The limit $\tu_{\tau_0}$ is also an embedding.
\item The convergence is in $C^{\infty}(S)$.
\end{itemize}
Another focus of this paper is studying how a family of embedded solutions can intersect with another solution or other pseudoholomorphic curves:\\
Assume $\tu_{\tau}=(a_{\tau},u_{\tau})$ , $-1<\tau\le 0$ is a continuous family of embedded pseudoholomorphic strips as in (\ref{main-boundary-value-problem}) with pairwise disjoint images. Let $\tilde{v}=(b,v)$ be either
\begin{enumerate}
\item another embedded solution of the boundary value problem (\ref{main-boundary-value-problem}) or
\item an embedded pseudoholomorphic disk with boundary condition $\{0\}\times {\mathcal D}^{\ast}$, i.e. a 'Bishop--disk' as in \cite{Hofer-Weinstein-conj} or \cite{AH},
\item a pseudoholomorphic half--cylinder over a periodic orbit $x(t)$ of $X_{\lambda}$ which lies on the surface ${\mathcal D}$, i.e. 
\[
\tilde{v}:Z^-:=(-\infty,0]\times S^1\longrightarrow {\bf R}\times M
\]
\[
\tilde{v}(s,t)=(s,x(t))\ ,\ x(S^1)\subset {\mathcal D}.
\]
\end{enumerate}
We study the intersection properties of the pseudoholomorphic curve $\tilde{v}$ with the family $\tu_{\tau}$. The following theorem states that 'there is no isolated first intersection'.

\begin{theorem}
\label{intersection-of-families-0}
Assume $\tu_{\tau}=(a_{\tau},u_{\tau})$ , $-1<\tau\le 0$ is a smooth family of embedded solutions of (\ref{main-boundary-value-problem}) with pairwise disjoint images. Let $\tilde{v}=(b,v)$ be as above. Moreover, we assume that $\tilde{u}_{\tau}$ and $\tilde{v}$ have disjoint images for $\tau<0$, but the intersection of $\tu_0({\bf R}\times [0,1])$ with $\tilde{v}({\bf R}\times [0,1])$ in case 1., $\tilde{v}(D)$ in case 2. and $\tilde{v}(Z^-)$ in case 3. is not empty. In the cases 2. and 3. the image of $\tilde{v}$ is contained in the image of $\tu_0$ or vice versa.\\
In case 1. this also holds unless the first intersection occurs at the boundary ${\bf R}\times\{0\}$, i.e. if $\tu_0(p)=\tilde{v}(q)$ for $p,q\in{\bf R}\times\{0\}$, and $\pas u_0(p)$, $\pas v(q)\in T_{u_0(p)}{\mathcal L}$ do not have the same orientation. In this case we can only conclude that $\tu_0({\bf R}\times \{0\})=\tilde{v}({\bf R}\times\{0\})$.
\end{theorem}

Since we will use the results and the notation from the papers \cite{part1} and \cite{part2}, we briefly summarize what is needed. As we mentioned earlier, we may modify the surface ${\mathcal D}$ near its boundary in order to achieve some normal form in local coordinates.
\begin{proposition}\label{3.1.1.}
Let $(M,\lambda)$ be a three-dimensional contact manifold. Further, let ${\mathcal L}$ be
a Legendrian knot and ${\mathcal D}$ an embedded surface with
$\partial{\mathcal D}={\mathcal L}$ so that all the singular points are non--degenerate. We denote the finitely many singular points on the boundary by $e_k$ , $1\le k\le N$ (ordered by orienting ${\mathcal L}$).\\ 
Then there is an embedded surface ${\mathcal D}'$
having the same boundary as ${\mathcal D}$ which differs from
${\mathcal D}$ only by a smooth $C^0$--small perturbation supported near
${\mathcal L}$ having the same singular points as ${\mathcal D}$ so that the following holds:\\
There is a neighborhood $U$ of ${\mathcal L}$ and a
diffeomorphism $\Phi:U\rightarrow S^1\times{\bf R}^2$ so that
\begin{itemize}
\item $\Phi^{\ast}(dy+xd\theta)=\lambda|_U\ ,\ (\theta,x,y)\in
S^1\times{\bf R}^2,$
\item $\Phi({\mathcal L})=S^1\times\{(0,0)\},$
\item $\Phi(e_k)=(\theta_k,0,0)$\ ,\ $0\le\theta_1<\cdots<\theta_N<1$,
\item $\Phi(U\cap{\mathcal
D}')=\{(\theta,a(\theta)r,b(\theta)r)\in S^1\times{\bf
R}^2\,|\,\theta,r\in[0,1]\},$
\end{itemize}
where $a,b$ are smooth 1--periodic functions with:
\begin{enumerate}
\item $b(\theta_k)=0$ and $b(\theta)$ is nonzero if $\theta\neq\theta_k$,
\item $a(\theta_k)<0$ if $e_k$ is a positive singular point, $a(\theta_k)>0$ if $e_k$ is a negative singular point,
\item if $e_k$ is elliptic then $-1<\frac{b'(\theta_k)}{a(\theta_k)}<0$,
\item if $e_k$ is hyperbolic then the quotient $\frac{b'(\theta_k)}{a(\theta_k)}$ is either strictly smaller than $-1$ or positive,
\item $a$ has exactly one zero in each of the intervals
$[\theta_k,\theta_{k+1}]$\,,\,$k=1,\ldots,N-1$ and $[\theta_N,1]\cup[0,\theta_1]$,
\item if $e_k$ is an elliptic singular point and if $|\theta-\theta_k|$ is sufficiently small then we have $b(\theta)=-\frac{1}{2}a(\theta)(\theta-\theta_k).$
\end{enumerate}
\end{proposition}
{\bf Proof:} See \cite{part1}.\qed

We showed in \cite{part1} that solutions $\tu$ to equation (\ref{main-boundary-value-problem}) with $\mbox{dist}(u({\bf R}\times\{1\}),\Gamma)>0$ converge to points $(0,p_{\pm})\in\{0\}\times{\mathcal L}\backslash\Gamma$ as $s\rightarrow\pm\infty$. It is sometimes convenient to modify the coordinates given by the above proposition near the points $p_{\pm}$ in order to make the surface ${\mathcal D}$ flat. Away from the boundary singular points $e_k$ we introduce the coordinate transformation
\begin{equation}\label{make-boundary-conditions-flat}
{\bf R}\times S^1\times {\bf R}^2\ni
(\tau,\theta,x,y)\longmapsto\left(\tau,\theta,x-\frac{a(\theta)}{b(\theta)}y,y\right)=(\tau,\theta,x-q(\theta)y,y).
\end{equation}
We then obtain the following coordinates on suitable neighborhoods $V_{\pm}$ of the points $p_{\pm}\in{\mathcal L}$:
\begin{equation}\label{make-boundary-conditions-flat-2}
\psi_{\pm}:{\bf R}^4\supset B_{\eps}(0)\widetilde{\longrightarrow}V_{\pm}\subset\RM,
\end{equation}
\[
\psi_{\pm}(0)=p_{\pm},
\]
\[
\psi_{\pm}({\bf R}^2\times\{0\}\times\{0\})=(\RL)\cap V_{\pm},
\]
\[
\psi_{\pm}(\{0\}\times{\bf R}\times\{0\}\times{\bf R}^{\pm})=(\{0\}\times {\mathcal D})\cap V_{\pm}.
\]
Using the coordinates $(\tau,\theta,x,y)$ for ${\bf R}^4$, the contact form on $\{0\}\times{\bf R}^3$ is then given by
\[
\hat{\lambda}_{\pm}=\psi_{\pm}^{\ast}\lambda=dy+\left(x+q(\theta)y\right)d\theta\ ,\ q(\theta):=\frac{a(\theta)}{b(\theta)}
\]
with Reeb vector field
\[
X_{\hat{\lambda}_{\pm}}=\frac{\partial}{\partial y}-q(\theta)\frac{\partial}{\partial x}
\]
(recall that the functions $a,b$ determine how the surface ${\mathcal D}$ is wrapping itself around the knot ${\mathcal L}$, see proposition \ref{3.1.1.}). Let $v_{\pm}(s,t):=(\psi^{-1}_{\pm}\circ \tu_0)(s,t)$ be the representative of a solution $\tu$ of (\ref{main-boundary-value-problem}) in the above coordinates for large $|s|$ (we also assume that $\mbox{dist}(u({\bf R}\times\{1\}),\Gamma)>0$). Our differential equation (\ref{main-boundary-value-problem}) has the following form in the above coordinates:
\[
v=(\tau,\theta,x,y):[s_0,\infty)\times[0,1]\To{\Real}^4
\]
\begin{equation}\label{5.1.5e.}
\partial_s v+M(v)\partial_t v=0
\end{equation}
\[v(s,0)  \in L_0={\Real}^2\times\{0\}\times\{0\} \]
\[ v(s,1)  \in L_1=\{0\}\times{\Real}\times\{0\}\times{\Real},\]
where $M$ is a suitable $4\times 4$--matrix valued function with $M^2=-\mbox{Id}$. We have shown in \cite{part1}:
\begin{theorem}\label{5.1.5a.}
There exist numbers $\rho,s'>0$ so that we have the following
estimate for each multi index $\alpha\in{\bf N}^2$ , $|\alpha|\ge
0$ and $s\ge s'$:
\[
\sup_{t\in [0,1]}|\partial^{\alpha}v(s,t)|\,\le
c_{\alpha}e^{-\rho(s-s')},
\]
where $c_{\alpha}$ are suitable positive constants.
\end{theorem}
\qed 
The main result of \cite{part1} is the following asymptotic formula for non constant solutions $v$ of (\ref{5.1.5e.}) having finite energy: 
\begin{theorem}\label{asymptotic-formula-theorem}
For sufficiently large $s_0$ and $|s|\ge s_0$ we have the following
asymptotic formulae for non constant solutions $v$ of
(\ref{5.1.5e.}) having finite energy:
\begin{equation}\label{asymptotic-formula}
v(s,t)=e^{\int_{s_0}^s\alpha_{\pm}(\tau)d\tau}\Big(e_{\pm}(t)+r_{\pm}(s,t)\Big),
\end{equation}
where $\alpha_{\pm}:[s_0,\infty)\rightarrow{\bf R}$ are smooth functions
satisfying $\alpha_+(s)\rightarrow\lambda_+ < 0$ and $\alpha_-(s)\rightarrow\lambda_->0$ as
$s\rightarrow\pm\infty$ with $\lambda_{\pm}\in {\bf Z}\frac{\pi}{2}$ being eigenvalues of the
selfadjoint operators
\[
A_{\pm\infty}:L^2([0,1],{\bf R}^4)\supset H^{1,2}_L([0,1],{\bf
R}^4)\longrightarrow L^2([0,1],{\bf R}^4)
\]
\[
\gamma\longmapsto -M_{\pm\infty}\dot{\gamma}\ ,\
M_{\pm\infty}:=\lim_{s\pm\rightarrow\infty}M(v(s,t)).
\]
Moreover, $e_{\pm}(t)$ is an eigenvector of $A_{\pm\infty}$ belonging to
the eigenvalue $\lambda_{\pm}$ with $e_{\pm}(t)\neq 0$ for all $t\in[0,1]$, 
and $r_{\pm}$ are smooth functions so that $r_{\pm}$
and all their derivatives converge to zero uniformly in $t$ as
$s\rightarrow\pm\infty$.
\end{theorem}
{\bf Proof:} See \cite{part1} \qed
The domain of the operators $A_{\pm\infty}$ above is the following dense subspace of $L^2([0,1],{\bf R}^4)$:
\[
H_L^{1,2}([0,1],{\bf R}^4):=\{\gamma\in H^{1,2}([0,1],{\bf
R}^4)\ |\ \gamma(0)\in L_0\,,\gamma(1)\in L_1\},
\] 
where 
\[
L_0:={\bf R}^2\times\{0\}\times\{0\}\ \mbox{and}\ L_1:=\{0\}\times{\bf R}\times\{0\}\times{\bf R}.
\]
In view of the Sobolev embedding theorem this definition makes sense. There is also the following refinement of the above asymptotic formula:
\begin{theorem}\label{convergence-of-alpha}
Let $v$ be as in theorem \ref{asymptotic-formula-theorem}. Then there is a constant $\delta>0$ such that for each integer $l\ge 0$ and each multi--index $\beta\in{\bf N}^2$
\[
\sup_{0\le t\le 1}|D^{\beta}r_{\pm}(s,t)|\ ,\ \left|\frac{d^l}{ds^l}(\alpha_{\pm}(s)-\lambda_{\pm})\right|\le c_{\beta,l} e^{-\delta|s|}
\]
with suitable constants $c_{\beta,l}>0$.
\end{theorem}
{\bf Proof:} See \cite{part1} \qed 

We recall from \cite{part1} that we chose the complex structure $J:\xi\rightarrow\xi$ near the Legendrian knot ${\mathcal L}$ as follows (in the local coordinates given by proposition \ref{3.1.1.}):
\begin{equation}\label{local-formula-for-J}
J(\theta,x,y)\cdot(1, 0, -x):= (0,-1,0)\ ,\ J(\theta,x,y)\cdot(0,1,0):=(1,0,-x).
\end{equation}
In the coordinates (\ref{make-boundary-conditions-flat}) the almost complex structure $\hat{J}$ on ${\bf R}^4$ induced by $\tilde{J}$ is given by 
\begin{equation}\label{matrix-of-J-hat}
\hat{J}(\tau,\theta,x,y) =  \left(\begin{array}{cc}
0 & -(x+q(\theta)y) \\ 0 & yq'(\theta) \\
-q(\theta) & -1+yq'(\theta)((x+q(\theta)y)q(\theta)-yq'(\theta)) \\1 & -(x+q(\theta)y)yq'(\theta)
\end{array}\right. 
\end{equation}
\[
\left.\begin{array}{cc}
0 & -1 \\1 & q(\theta) \\(x+q(\theta)y)q(\theta)-yq'(\theta) & q(\theta)((x+q(\theta)y)q(\theta)-yq'(\theta)) \\
-(x+q(\theta)y) & -(x+q(\theta)y)q(\theta)
\end{array}\right).
\]
If $\lambda_{\pm}$ is an odd integer multiple of $\pi/2$ then the asymptotic formula of theorem \ref{asymptotic-formula-theorem} looks as follows:
\begin{eqnarray}\label{explicit-asymptotic-formula}
v_{\pm}(s,t) & = & -{\kappa}_{\pm}e^{\int_{s_0}^s\alpha_{\pm}(\tau)d\tau}\Big(\cos(\lambda_{\pm} t), -q_{\pm}(0)\cos(\lambda_{\pm} t), 0, \sin(\lambda_{\pm} t)\Big)+ \nonumber \\
 & & +e^{\int_{s_0}^s\alpha_{\pm}(\tau)d\tau}\ve_{\pm}(s,t).
\end{eqnarray}
In the following we will denote by $\varepsilon(s,t)$ any ${\bf R}^4$-- or real--valued function which converges to zero with all its derivatives uniformly in $t$ as $s\rightarrow\pm\infty$ if we are not interested in the particular function. In order to simplify notation we will often drop the subscript $\pm$. Using the fact that $\alpha'(s)\rightarrow 0$ as $|s|\rightarrow \infty$ (proved in \cite{part1}), we obtain the following asymptotic formulae for the derivatives of $v(s,t)$ ($\kappa,\kappa_{\pm}$ are suitable nonzero constants):
\begin{eqnarray}\label{explicit-asymptotic-formula-for-us}
\partial_s v(s,t) & = & e^{\int_{s_0}^s\alpha(\tau)d\tau}
 \cdot\Big[-{\kappa}\Big(\lambda\cos(\lambda t), -\lambda q(0)\cos(\lambda t), 0, \lambda\sin(\lambda t)\Big)+\nonumber\\
  & & +\varepsilon(s,t)\Big],
\end{eqnarray}
\begin{eqnarray}\label{explicit-asymptotic-formula-for-ut}
\partial_t v(s,t) & = & e^{\int_{s_0}^s\alpha(\tau)d\tau}\cdot
 \Big[-{\kappa}\Big(-\lambda\sin(\lambda t), \lambda q(0)\sin(\lambda t), 0, \lambda\cos(\lambda t)\Big)+\nonumber\\
  & & +\varepsilon(s,t)\Big].
\end{eqnarray}

If we use the coordinates given by proposition \ref{3.1.1.} without making the boundary conditions 'flat' as in (\ref{make-boundary-conditions-flat}) then the appropriate versions of (\ref{explicit-asymptotic-formula}) and (\ref{explicit-asymptotic-formula-for-us}) are the following. If $\lambda_{\pm}$ is an odd integer multiple of $\pi/2$ we have:
\begin{eqnarray}\label{odd-explicit-asymptotic-formula-v2}
 & & v_{\pm}(s,t)\nonumber \\ 
& = & -{\kappa}_{\pm}e^{\int_{s_0}^s\alpha_{\pm}(\tau)d\tau}\Big(\cos(\lambda_{\pm} t), -q_{\pm}(0)\cos(\lambda_{\pm} t), q_{\pm}(0)\sin(\lambda_{\pm} t), \sin(\lambda_{\pm} t)\Big)+ \nonumber \\
 & & +e^{\int_{s_0}^s\alpha_{\pm}(\tau)d\tau}\ve_{\pm}(s,t)
\end{eqnarray}
and
\begin{eqnarray}\label{odd-explicit-asymptotic-formula-for-us-v2}
& & \partial_s v_{\pm}(s,t)\nonumber \\
 & = & e^{\int_{s_0}^s\alpha_{\pm}(\tau)d\tau}\cdot \Big[-{\kappa}_{\pm}\Big(\lambda_{\pm}\cos(\lambda_{\pm} t), -\lambda_{\pm} q_{\pm}(0)\cos(\lambda_{\pm} t),\\
 & &  \lambda_{\pm} q_{\pm}(0)\sin(\lambda_{\pm} t), \lambda_{\pm}\sin(\lambda_{\pm} t)\Big)+\varepsilon_{\pm}(s,t)\Big].\nonumber
\end{eqnarray}
For $\lambda_{\pm}\in{\bf Z}\pi$ we have
\begin{eqnarray}\label{even-explicit-asymptotic-formula-v2}
 & & v_{\pm}(s,t)\nonumber \\ 
& = & {\kappa}_{\pm}e^{\int_{s_0}^s\alpha_{\pm}(\tau)d\tau}\Big(0 , \cos(\lambda_{\pm} t), -\sin(\lambda_{\pm} t), 0\Big)+ \nonumber \\
 & & +e^{\int_{s_0}^s\alpha_{\pm}(\tau)d\tau}\ve_{\pm}(s,t)
\end{eqnarray}
and
\begin{eqnarray}\label{even-explicit-asymptotic-formula-for-us-v2}
\partial_s v_{\pm}(s,t) & = & e^{\int_{s_0}^s\alpha_{\pm}(\tau)d\tau}\cdot \nonumber\\
 & & \cdot\Big[{\kappa}_{\pm}\Big(0 , \lambda_{\pm} \cos(\lambda_{\pm} t), -\lambda_{\pm}\sin(\lambda_{\pm} t),0\Big)+\varepsilon_{\pm}(s,t)\Big].
\end{eqnarray}
The following theorem is the main result of the paper \cite{part2}:
\begin{theorem}\label{main-implicit-function-theorem}
{\bf (Implicit function theorem)}\\
Let $\tu_0=(a_0,u_0)$ be an embedded solution of (\ref{main-boundary-value-problem}) so that its Maslov--index $\mu(\tu_0)$ vanishes and $\mbox{dist}(u_0({\bf R}\times\{1\}),\Gamma)>0$. Assume moreover, that $|\tu_0(s,t)-p_{\pm}|$ decays either like $e^{-\pi|s|}$ or like $e^{-\frac{\pi}{2}|s|}$ for large $|s|$ in local coordinates near the points $p_{\pm}:=\lim_{s\rightarrow\pm\infty}\tu_0(s,t)$ and that $p_-\neq p_+$. Then there is a smooth family $(\tilde{v}_{\tau})_{-1<\tau<1}$ of embedded solutions of (\ref{main-boundary-value-problem}) with the following properties:
\begin{itemize}
\item $\tilde{v}_0=\tu_0$,
\item The solutions $\tilde{v}_{\tau}$ have the same Maslov--index and the same decay rates as $\tu_0$,
\item The sets 
\[
U_{\pm}:=\bigcup_{-1<\tau< 1}\{\lim_{s\rightarrow\pm\infty}\tilde{v}_{\tau}(s,t)\}
\]
are open neighborhoods of the points $p_{\pm}$ in ${\mathcal L}$.
\end{itemize}
If $|\tu_0(s,t)-p_{\pm}|$ decays like $e^{-\frac{\pi}{2}|s|}$ for both $s\rightarrow+\infty$ and $s\rightarrow-\infty$ then we have in addition
\begin{itemize}
\item $\tilde{v}_{\tau}(S)\cap\tilde{v}_{\tau'}(S)=\emptyset$ if $\tau\neq\tau'$.
\end{itemize}
\end{theorem}
{\bf Proof:} See \cite{part2}
\qed
 
In section \ref{IFT-second-version}, we will prove the following version of theorem \ref{main-implicit-function-theorem}:
\begin{theorem}
{\bf (Implicit Function Theorem--second version)}\\
Let $\tu_0=(a_0,u_0)$ be an immersed solution of (\ref{main-boundary-value-problem}). Denote by $\lambda_{\pm}$ the decay rates of $\tu_0$ as in theorem \ref{asymptotic-formula-theorem} and let $\mu(\tu_0)$ be the Maslov index of $\tu_0$. Assume that one of the following conditions are satisfied:
\begin{enumerate}
\item $\lambda_{\pm}=\mp m_{\pm}\pi$ with integers $m_{\pm}\ge 1$ and $\mu(\tu_0)<-\frac{1}{2}(m_-+m_+)$ holds,
\item the absolute value of one of the numbers $\lambda_{\pm}$ equals $\frac{\pi}{2}$, and the absolute value of the other equals $m\pi$ with some positive integer $m$. We also assume that $\mu(\tu_0)<-\frac{1}{2}m-\frac{1}{4}$.
\end{enumerate}
Assume moreover that $\mbox{dist}(u_0({\bf R}\times\{1\}),\Gamma)>0$. Then there is an integer $N\ge 1$ and a smooth family $(\tilde{v}_{\tau})_{\tau\in {\bf R}^N}$ of solutions the boundary value problem (\ref{main-boundary-value-problem}) with the following properties:
\begin{itemize}
\item $\tilde{v}_0=\tu_0$ and $\tilde{v}_{\tau}\not\equiv\tu_0$ if $\tau\neq 0$
\item the solutions $\tilde{v}_{\tau}$ have the same end points as the solution $\tu_0$, i.e.
\[
\lim_{s\rightarrow\pm\infty}\tu_0(s,t)=\lim_{s\rightarrow\pm\infty}\tilde{v}_{\tau}(s,t)\ \forall\ \tau,
\]
\item The solutions $\tilde{v}_{\tau}$ have the same Maslov--index and the same decay rates as $\tu_0$,
\end{itemize}
\end{theorem}

We can see now the purpose of theorem \ref{compactness-result}: The 'implicit function theorems' can only be applied to embedded or immersed solutions $\tu_0$ where Maslov--index and decay rates are related in a certain way. We saw in the paper \cite{part2} that otherwise the underlying Fredholm operator would have negative index. Our compactness result, theorem \ref{compactness-result} has to make sure that theorems \ref{main-implicit-function-theorem} and \ref{extended-implicit-function-theorem} remain applicable to the $C^{\infty}_{loc}$--limit of a sequence of solutions. We have proved also the following result in \cite{part2} which is based on the  maximum principle:
\begin{proposition}\label{maximum-principle}
Let $\tu=(a,u):S\longrightarrow\RM$ be a non--constant solution of
the boundary value problem (\ref{main-boundary-value-problem}). 
\begin{itemize}
\item Then the path
$s\longmapsto u(s,1)$ is transverse to the characteristic
foliation, i.e. $\pas u(s,1)\not\in\ker\lambda(u(s,1))$. We actually have
\[
0<\pat a(s,1)=-\lambda(u(s,1))\pas u(s,1)
\]
for all $s\in{\bf R}$.
\item We have $a(s,t)<0$ whenever $0\le t<1$,
\item The pseudoholomorphic strip never hits $\{0\}\times{\mathcal L}$, i.e.
\[
\tu(S)\cap (\{0\}\times{\mathcal L})=\emptyset.
\]
In particular,
\[
\lim_{s\rightarrow\pm\infty}\tu(s,t)\not\in\tu(S).
\]
\end{itemize}
\end{proposition}
{\bf Proof:}
See \cite{part2}.\qed

\section{Local existence of solutions}\label{local-existence-of-sol}

In this section we will prove theorem \ref{local-existence-theorem}, i.e. we establish local fillings by pseudoholomorphic curves near an elliptic singularity at the boundary.
 Because of proposition \ref{3.1.1.} we are in the following situation near an elliptic singular point $e\in{\mathcal L}$:\\
We may assume that the contact manifold is the three dimensional
Euclidean space $\{(\theta,x,y)\in{\bf R}^3\}$ endowed with the
contact form $\lambda=dy+xd\theta$. The piece of the Legendrian
knot situated near $e$ corresponds to some interval
$\{(\theta,0,0)\in{\bf R}^3\,|\,|\theta|<\varepsilon\}$, where
$\varepsilon>0$ is a suitable constant. The elliptic
singular point then corresponds to the origin in ${\bf R}^3$ and the spanning surface
${\mathcal D}$ is given by $\{(\theta,x,-\frac{1}{2}\theta
x)\in{\bf R}^3\,|\,|\theta|<\varepsilon\,,\,x\le 0\}$ if $e$ is a positive elliptic point, otherwise we have $x\ge 0$. We start constructing solutions near $e$. The contact structure is generated by the
vectors
\[e_1=\left(\begin{array}{c} 1 \\ 0 \\ -x \end{array}\right)\ \mbox{and}\ e_2=\left(\begin{array}{c} 0 \\ 1 \\ 0 \end{array}\right).\]
We have chosen a particular complex structure $J$ on $\ker\lambda$ near the Legendrian knot by demanding
\begin{equation}\label{standard-J-near-e+-}
Je_1:=-e_2\ \mbox{and}\ Je_2=e_1.
\end{equation}
This complex structure is compatible with $d\lambda$, i.e.
$d\lambda\circ(\mbox{Id}\times J)$ is a bundle metric and defines an almost
complex structure $\tilde{J}$ in the usual way.\\ The
boundary value problem, we are going to study, is the following:
\[
 \begin{array}{ll}
          \tilde{u}=(a,u):S\longrightarrow\RM & \\
          \partial_s\tu+\tilde{J}(\tu)\partial_t\tu=0 & \\
          \tu(s,0)\subset{\bf R}\times{\cal L} & \\
          \tu(s,1)\subset\{0\}\times{\cal D}^{\ast} & \\
           0< E(\tu)<+\infty & 
          \end{array},
\]
where ${\mathcal D}^{\ast}$ is
the spanning surface without the singular points, $S:={\bf
R}\times[0,1]$ and
\[E(\tu)=\sup_{\phi\in\Sigma}\int_S\tu^{\ast}d(\phi\lambda)\
\mbox{ (energy of $\tu$ ) }\] with $\Sigma:=\{\phi\in
C^{\infty}({\bf R},[0,1])\,|\,\phi'\ge 0\}$. Since we have chosen
good coordinates near the elliptic singular points and an explicit
almost complex structure $\tilde{J}$, we will be able to explicitly state 
1--parameter families of solutions to the above boundary value
problem near the elliptic singular points.\\ These solutions look simpler after having performed a biholomorphic transformation of the domain as
follows: Let $\Omega:=\{z=s+it\in{\bf C}\,|\,s^2+t^2\le 1\,,\,t\ge
0\}\backslash\{-1,+1\}$ be the upper half disk in the complex
plane without the corner points. The infinite strip $S$ and
$\Omega$ are equivalent via the biholomorphic map
\[S\longrightarrow\Omega\]
\begin{equation}\label{strip-halfdisk}
s+it\longmapsto\frac{e^{\frac{\pi}{2}(s+it)}-1}{e^{\frac{\pi}{2}(s+it)}+1}=\mbox{tanh}\left(\frac{\pi}{4}(s+it)\right).
\end{equation}
Under this transformation, ${\bf R}\times\{0\}$ is mapped onto
$(-1,+1)$ and ${\bf R}\times\{1\}$ is mapped onto
$\{s+i\sqrt{1-s^2}\in{\bf C}\,|\,s\in (-1,+1)\}$. We write in
coordinates
\[\tilde{u}=(a,\theta,x,y):\Omega\rightarrow{\bf R}\times{\bf R}^3\]
and obtain the following boundary value problem:
\begin{eqnarray*}
\pas a-\pat y-x\pat \theta & = & 0 \\
\pas\theta+\pat x & = & 0 \\
\pas x-\pat\theta & = & 0 \\
\pat a+\pas y+x\pas \theta & = & 0 \\
x(s,0)\equiv y(s,0) & \equiv & 0 \\
a(s,\sqrt{1-s^2}) & \equiv & 0 \\
y(s,\sqrt{1-s^2}) & = & -\frac{1}{2}(x\theta)(s,\sqrt{1-s^2}).
\end{eqnarray*}
The following maps satisfy the above boundary
value problem as long as they stay in the coordinate patch near
the elliptic singular point where the Seifert surface and $J$ are in normal form:
\begin{equation}\label{local_solutions}
\tu_{\varepsilon}(s,t)=\left(\frac{1}{4}\varepsilon^2(s^2+t^2-1),\varepsilon
s,-\varepsilon t,\frac{1}{2}\varepsilon^2st\right),
\end{equation}
with $\ve>0$ if $e$ is a positive elliptic singular point and $\ve<0$ otherwise. 
Transforming back the infinite strip $S={\bf R}\times[0,1]$ the solutions (\ref{local_solutions}) become
\begin{eqnarray}\label{local-solutions-2}
\tu_{\ve}(s,t) & = & \left(
-\frac{\ve^2\cos\left(\frac{\pi t}{2}\right)}{2\left[\cos\left(\frac{\pi t}{2}\right)+\cosh\left(\frac{\pi s}{2}\right)\right]}, \frac{\ve \sinh\left(\frac{\pi s}{2}\right)}{\cos\left(\frac{\pi t}{2}\right)+\cosh\left(\frac{\pi s}{2}\right)}, \right. \nonumber \\
 & & \left.\frac{-\ve\sin\left(\frac{\pi t}{2}\right)}{\cos\left(\frac{\pi t}{2}\right)+\cosh\left(\frac{\pi s}{2}\right)}, \frac{\ve^2\sin\left(\frac{\pi t}{2}\right)\sinh\left(\frac{\pi s}{2}\right)}{2\left[\cos\left(\frac{\pi t}{2}\right)+\cosh\left(\frac{\pi s}{2}\right)\right]^2}
\right).
\end{eqnarray}
These solutions obviously satisfy all the requirements of theorem \ref{local-existence-theorem}, and they decay like $e^{-\frac{\pi}{2}|s|}$ near the ends. By the asymptotic formula, theorem \ref{asymptotic-formula-theorem}, this is the slowest possible decay rate.

\section{Local reflection of solutions at the boundary}

The following lemma provides convenient local coordinates near a boundary point of a solution. We obtain as a corollary that pseudoholomorphic curves with totally real boundary conditions can be locally 'Schwarz--reflected' near a regular boundary point. In the special case where the boundary condition is ${\bf R}\times {\mathcal L}\subset {\bf R}\times M$ we can locally reflect in any boundary point. We confine ourselves to the case of dimension four for notational simplicity.
\begin{lemma}\label{local-coordinates-near-boundary-point}
Let $(W,J)$ be an almost complex manifold of dimension four and let $F$ be a totally real submanifold. Furthermore, assume that $u_0:D^+\rightarrow W$ is an embedded $J$--holomorphic half--disk with boundary condition $u_0(D^+\cap{\bf R})\subset F$, where $D^+:=\{z\in{\bf C}\,|\,\mbox{Im}(z)\ge 0\ ,\ |z|\le 1\}$. Then there are $0<\ve\le 1$, a neighborhood $U$ of $p=u_0(0)$ in $W$ and a diffeomorphism $\phi:U\rightarrow V\subset{\bf C}^2$ onto an open neighborhood $V$ of $0$ in ${\bf C}^2$ such that
\begin{enumerate}
\item $\phi(U\cap F)=V\cap {\bf R}^2$,
\item $(\phi\circ u_0)(z)=(z,0)\ ,\ z\in(D^+_{\ve}\times \{0\})\cap V$, where $D^+_{\ve}:=\{z\in{\bf C}\,|\,\mbox{Im}(z)\ge 0\ ,\ |z|\le \ve\}$.
\item The induced almost complex structure $\bar{J}(q)=D\phi(\phi^{-1}(q))\circ J(\phi^{-1}(q))\circ D\phi^{-1}(q)$ satisfies
\[
\bar{J}(z,0)=i\ \mbox{if}\ (z,0)\in({\bf C}\times\{0\})\cap V.
\]
\end{enumerate}
\end{lemma}

{\bf Proof:}\\
Since $u$ is embedded, there is a coordinate map $\phi$ which satisfies condition 2. Hence we may assume that $W={\bf C}^2$ and $u_0(z)=(z,0)$. We denote $\phi(U\cap F)$ again by $F$. Points $(s,0)\in{\bf R}\times \{0\}\subset {\bf C}^2$ near $0$ are then contained in $F$. The map $u_0(z)=(z,0)$ is $\bar{J}$--holomorphic. If we pick
\[
(k,0)=Du_0(z)k\in {\bf C}\times\{0\}=T_{u_0(z)}(u_0(D^+))
\]
then
\begin{eqnarray*}
\bar{J}(u_0(z))(k,0) & = & \bar{J}(u_0(z))Du_0(z)k \\
 & = & Du_0(z)(ik) \\
 & = & i(k,0),
\end{eqnarray*}
hence $\bar{J}(u_0(z))$ acts as multiplication by $i$ on ${\bf C}\times \{0\}$. We will successively change coordinates until the other two conditions are satisfied as well. Let us take care of condition 3. first. We can find a smooth map $\psi:{\bf C}\rightarrow \mbox{GL}_{{\bf R}}({\bf C}^2)$ into the set of real linear automorphisms of ${\bf C}^2$ with the properties:
\begin{enumerate}
\item $\psi(z)\bar{J}(u_0(z))=i\,\psi(z)$ for all $z\in {\bf C}$,
\item $\psi(z)(h,0)=(h,0)$ for all $h\in {\bf C}$ and $z\in D^+$.
\end{enumerate}
If $s\in {\bf R}$ then $T_{u_0(s)}F=T_{(s,0)}F$ is generated by the vector $(1,0)$ and some vector valued function $(i\alpha,\beta)(s)$ for suitable $\alpha(s)\in{\bf R}$ and $0\neq \beta(s)\in{\bf C}$. We may demand in addition to the above two conditions that $\psi(s)(i\alpha(s),\beta(s))=(0,1)$ for real $s$ so that
\[
\psi(s)T_{u_0(s)}F={\bf R}^2\subset{\bf C}^2.
\]
Define now a map $\sigma:V'\rightarrow V''$ between suitable neighborhoods $V',V''$ of $(0,0)$ in ${\bf C}^2$ by
\[
\sigma(v,w):=\psi(v)(v,w).
\]
We have $\sigma(z,0)=(z,0)$ for all $z\in D^+$ and also  $D\sigma(z,0)(h,k)=\psi(z)(h,k)$ for $z\in D^+$. We conclude that $\sigma$ is a diffeomorphism onto its image, provided $V',V''$ are chosen sufficiently small. We then restrict $u_0$ to a smaller half--disk $D^+_{\ve}$ so that $u_0(D^+_{\ve})\subset V'$. We have also arranged that 
\[
D\sigma(s,0)(T_{u_0(s)}F)=\psi(s)(T_{u_0(s)}F)={\bf R}^2
\]
for $s\in{\bf R}$. Composing the original coordinate map $\phi$ with $\sigma$ we may assume now that conditions 2. and 3. of the lemma are satisfied. Moreover, we have 
\[
T_{u_0(s)}F={\bf R}^2 \ \mbox{for}\  s\in{\bf R}\cap D^+. 
\]
We will finally achieve condition 1. by another modification. We will find a local diffeomorphism $\chi$ with the properties
\begin{itemize}
\item $\chi(z,0)=(z,0)$,
\item $\chi(F)\subset {\bf R}^2$, at least for a piece of $F$ near ${\bf R}\times\{0\}$,
\item $D\chi(z,0)$ is complex linear.
\end{itemize}
The first and the last items above ensure that we maintain properties 2. and 3. Since ${\bf R}\times\{0\}\subset F$ and $T_{(s,0)}F={\bf R}^2$ we may write $F$ near ${\bf R}\times\{0\}$ as
\[
F=\{(s+if(s,t),t+ig(s,t))\in{\bf C}^2\,|\,s,t\in{\bf R}\,\}
\]
for suitable real valued smooth functions $f,g$ which satisfy in addition
\[
f(s,0)\equiv g(s,0)\equiv 0\ ,\ \partial_tf(s,0)\equiv\partial_tg(s,0)\equiv 0.
\]
We define $f,g$ also for complex arguments by $f(z,w):=f(\mbox{Re}(z),\mbox{Re}(w))$ and similarly $g$. We define the map $\chi$ as follows:
\[
\chi(z,w):=(z-if(z,w),w-ig(z,w)).
\]
This map has the desired properties.\qed

We note the following corollary:
\begin{corollary}\label{local-reflection}
Let $(W,J)$ be an almost complex manifold of dimension four and let $F$ be a totally real submanifold. Furthermore, assume that $u_0:D^+\rightarrow W$ is a $J$--holomorphic half--disk with boundary condition $u_0(D^+\cap{\bf R})\subset F$ so that $\pas u_0(0)\neq 0$. Then there are $0<\ve$ and a $J$--holomorphic disk $v_0:D_{\ve}\rightarrow W$ defined on the full disk of radius $\ve$ so that $v_0|_{D^+_{\ve}}\equiv u_0|_{D^+_{\ve}}$.
\end{corollary}
{\bf Proof:}\\
We use lemma \ref{local-coordinates-near-boundary-point}. Then the map $z\mapsto (z,0)$ is $J$--holomorphic on the full disk because $J(z,0)\equiv i$ for all $z\in{\bf C}$.
\qed
We return to the special situation where $(W,J)=({\bf R}\times M,\tilde{J})$ where $\tilde{J}$ has the form (\ref{local-formula-for-J}) near the Legendrian knot ${\mathcal L}$, and where the totally real submanifold $F$ is ${\bf R}\times {\mathcal L}$.
\begin{proposition}\label{reflection-2}
For the special choices of $W,\tilde{J}$ and $F$ outlined above the assertion of corollary \ref{local-reflection} also holds without the assumption $\pas u_0(0)\neq 0$.
\end{proposition}
{\bf Proof:}\\
In local coordinates the map $u_0=(a,\theta,x,y):D^+\rightarrow {\bf R}^4$ satisfies the following differential equation (see section \ref{local-existence-of-sol}):
\begin{eqnarray*}
\pas a-\pat y-x\pat \theta & = & 0 \\
\pas\theta+\pat x & = & 0 \\
\pas x-\pat\theta & = & 0 \\
\pat a+\pas y+x\pas \theta & = & 0 \\
x(s,0)\equiv y(s,0) & \equiv & 0.
\end{eqnarray*}
We then use the obvious Schwarz reflection.\qed

\section{Intersection Properties}

The following intersection result is more general than needed for this paper, but it will be useful later on. Assume $\tu_{\tau}=(a_{\tau},u_{\tau})$ , $-1<\tau\le 0$ is a smooth family of embedded pseudoholomorphic strips as in (\ref{main-boundary-value-problem}) with pairwise disjoint images. Let $\tilde{v}=(b,v)$ be either
\begin{enumerate}
\item another embedded solution of the boundary value problem (\ref{main-boundary-value-problem}) or
\item an embedded pseudoholomorphic disk with boundary condition $\{0\}\times {\mathcal D}^{\ast}$, i.e. a 'Bishop--disk' as in \cite{Hofer-Weinstein-conj} or \cite{AH},
\item a pseudoholomorphic half--cylinder over a periodic orbit $x(t)$ of $X_{\lambda}$ which lies on the surface ${\mathcal D}$, i.e. 
\[
\tilde{v}:Z^-:=(-\infty,0]\times S^1\longrightarrow {\bf R}\times M
\]
\[
\tilde{v}(s,t)=(s,x(t))\ ,\ x(S^1)\subset {\mathcal D}.
\]
\end{enumerate}
We study the intersection properties of the pseudoholomorphic curve $\tilde{v}$ with the family $\tu_{\tau}$. The following theorem states that 'there is no isolated first intersection'. 

\begin{theorem}\label{interior-and-boundary-intersections}
Assume $\tu_{\tau}=(a_{\tau},u_{\tau})$ , $-1<\tau\le 0$ is a smooth family of embedded solutions of (\ref{main-boundary-value-problem}) with pairwise disjoint images. Let $\tilde{v}=(b,v)$ be as above. Moreover, we assume that $\tilde{u}_{\tau}$ and $\tilde{v}$ have disjoint images for $\tau<0$, but the intersection of $\tu_0({\bf R}\times [0,1])$ with $\tilde{v}({\bf R}\times [0,1])$ in case 1., $\tilde{v}(D)$ in case 2. and $\tilde{v}(Z^-)$ in case 3. is not empty. In the cases 2. and 3. the image of $\tilde{v}$ is contained in the image of $\tu_0$ or vice versa.\\
In case 1. this also holds unless the first intersection occurs at the boundary ${\bf R}\times\{0\}$, i.e. if $\tu_0(p)=\tilde{v}(q)$ for $p,q\in{\bf R}\times\{0\}$, and $\pas u_0(p)$, $\pas v(q)\in T_{u_0(p)}{\mathcal L}$ do not have the same orientation. Then we can only conclude that $\tu_0({\bf R}\times \{0\})=\tilde{v}({\bf R}\times\{0\})$.
\end{theorem}

The conclusion of the theorem is absurd in the cases 2. and 3.: In case 2. the conclusion of the theorem would violate proposition \ref{maximum-principle}. In case 3. it would imply that $\pil\pas u_0\equiv 0$ on an open subset of the domain. This implies $\pil\pas u_0\equiv 0$ on all of $S$ because the components of $\pil\pas u_0$ with respect to a complex frame for the bundle $(u_0^{\ast}\xi,J(u_0))$ satisfy a Cauchy Riemann type equation where the Similarity Principle (see appendix \ref{similarity-principle}) is applicable. In the paper \cite{part1} (lemma 4.1) we have shown that this implies that $\tu_0$ is constant in contradiction to the assumption that $\tu_0$ is an embedding. Hence we obtain the following corollary from the first part of theorem \ref{interior-and-boundary-intersections}:
\begin{corollary}\label{disks-cylinders}
Assume $\tu_{\tau}=(a_{\tau},u_{\tau})$ , $-1<\tau\le 0$ is a continuous family of embedded solutions of (\ref{main-boundary-value-problem}) with pairwise disjoint images. Let $\Gamma\subset{\mathcal D}$ be a curve which is either the trace $v(\partial D)$ of a Bishop--disk or periodic orbit $x(t)=\tilde{v}(0,t)$ of the Reeb vector field as above. Moreover, we assume that $\tilde{u}_{\tau}$ and $\tilde{v}$ do not intersect for $\tau<0$. Then the curve $u_0(s,1)$ on ${\mathcal D}$ and $\Gamma$ do not intersect either. 
\end{corollary}
\qed

The following two propositions are local versions of theorem \ref{interior-and-boundary-intersections}. Their proofs are very similar. Since the boundary case is more difficult we will only prove proposition \ref{intersection-of-families-local}.
\begin{proposition}\label{intersection-of-families-local-interior}
Let $u_{\tau}:D\rightarrow {\bf C}^2$ , $-1<\tau\le 0$ be a smooth family of embedded solutions of the differential equation
\begin{equation}
\partial_s u_{\tau}+J(u_{\tau})\partial_t u_{\tau}=0,
\end{equation}
where $J$ is an almost complex structure on ${\bf C}^2$ satisfying $J(z,0)=i$ and $D$ is the open unit disk in the complex plane. Assume furthermore that $u_0(z)=(z,0)$ and that the images of the $u_{\tau}$ are pairwise disjoint. Let $v$ be another embedded solution of the same differential equation with $u_{\tau}(D)\cap v(D)=\emptyset$ for $\tau<0$, but $v(0)=u_0(0)=0$. Then either $v(D)\subset u_0(D)=D\times\{0\}$ or $D\times\{0\}\subset v(D)$.
\end{proposition}
\qed

\begin{proposition}\label{intersection-of-families-local}
Let $u_{\tau}:D^+\rightarrow {\bf C}^2$ , $-1<\tau\le 0$ be a smooth family of embedded solutions of the boundary value problem
\begin{equation}\label{local-boundary-value-problem}
\partial_s u_{\tau}+J(u_{\tau})\partial_t u_{\tau}=0
\end{equation}
\[
u_{\tau}(D^+\cap{\bf R})\subset{\bf R}^2,
\]
where $J$ is an almost complex structure on ${\bf C}^2$ satisfying $J(z,0)=i$. Assume furthermore that $u_0(z)=(z,0)$ and that the images of the $u_{\tau}$ are pairwise disjoint. Let $v$ be another embedded solution of the boundary value problem (\ref{local-boundary-value-problem}) with $u_{\tau}(D^+)\cap v(D^+)=\emptyset$ for $\tau<0$, but $v(0)=u_0(0)=0$.\\ 
Denoting the upper/lower half--planes in ${\bf C}$ by $H^{\pm}$ and the upper half--disk with radius $r$ by $D^+_r$, we assume moreover that there is some $\ve>0$ such that $v(D^+_{\ve})\subset H^+\times{\bf C}$. Then either $v(D^+)\subset u_0(D^+)=D^+\times\{0\}$ or $D^+\times\{0\}\subset v(D^+)$. If $v(D^+_{\ve})\subset H^-\times{\bf C}$ then either $v(D^+)\subset D^-\times\{0\}$ or $D^-\times\{0\}\subset v(D^+)$, i.e. $v$ agrees with the Schwarz--reflection of $u_0$.
\end{proposition}

\begin{figure}[!ht]
\begin{center}
\epsfig{file=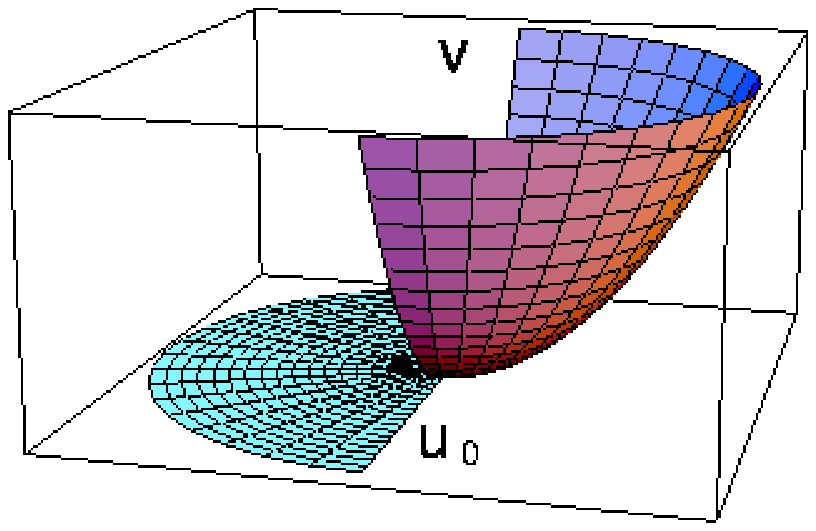,height=6cm,angle=0}
\caption{}\label{bad}
\end{center}
\end{figure}

\begin{figure}[!ht]
\begin{center}
\epsfig{file=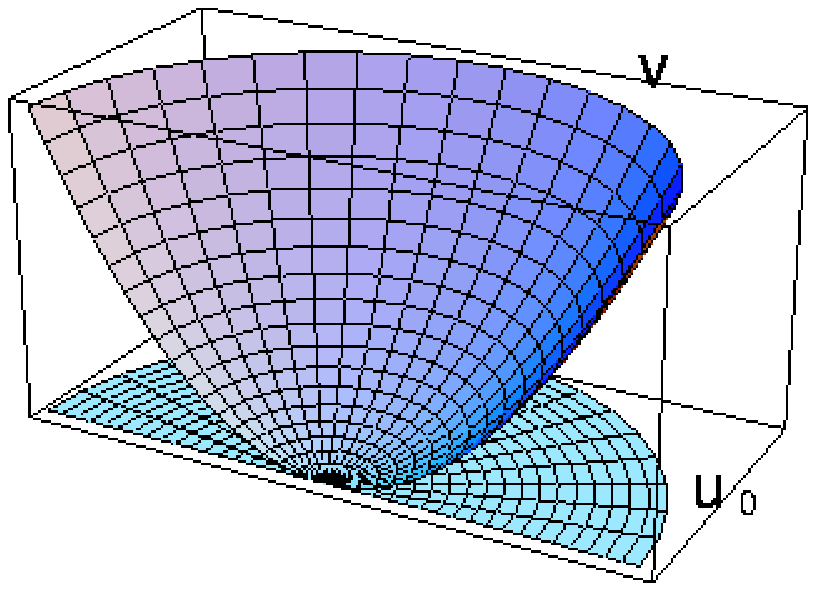,height=6cm,angle=0}
\caption{}\label{good}
\end{center}
\end{figure}

{\bf Remark:} The assumption that $v(D^+_{\ve})\subset H^+\times{\bf C}$ avoids the hypothetical intersection picture in figure \ref{bad}, where $v(D_{\ve}^+)$ lies on the 'wrong side'. The situation we are looking at in this case is the one in figure \ref{good}.

\begin{figure}[!ht]
\begin{center}
\epsfig{file=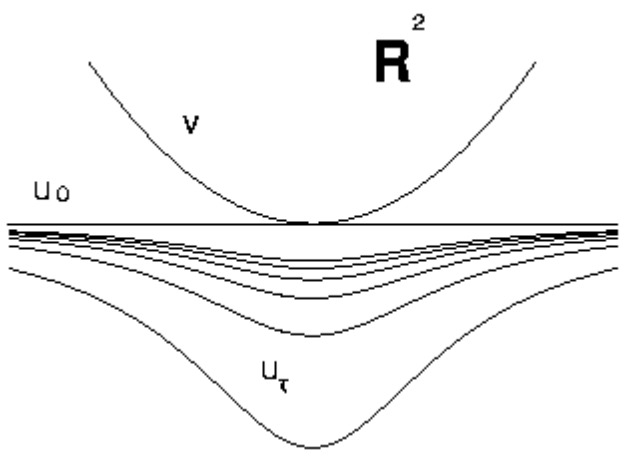,height=6cm,angle=0}
\caption{The graphs of $v|_{{\bf R}\cap D^+}$ and $u_{\tau}|_{{\bf R}\cap D^+}$.}\label{boundary-curves}
\end{center}
\end{figure}

{\bf Proof of proposition \ref{intersection-of-families-local}:}\\
Writing $v(z)=(a(z),b(z))$ it suffices to show that $b(z)$ vanishes on some neighborhood of $0$. We first note that $v$ and $u_0$ are tangent at $0$. Indeed, we have 
\[
\frac{d}{ds}v(s,0)|_{s=0}=Dv(0)\,1\in {\bf R}\times\{0\},
\]
otherwise the curve $v(s,0)$ would intersect some of the curves $u_{\tau}(s,0)$ for $\tau<0$ as well (see figure \ref{boundary-curves}), which does not happen by assumption. 
Then $Dv(0)\,i=i\,Dv(0)\,1\in i{\bf R}\times\{0\}$, hence the range of $Dv(0)$ is ${\bf C}\times\{0\}$. We have
\begin{equation}\label{taylor-expansion-of-v}
v(z)=Dv(0)z+O(|z|^2)=(c\,z,0)+O(|z|^2)
\end{equation}
for a suitable number $c\in{\bf C}\backslash\{0\}$, and even $c\in{\bf R}\backslash\{0\}$ because of $Dv(0)\,1\in {\bf R}\times\{0\}$. After maybe replacing $D^+$ by a smaller half--disk we may assume that the image of $v$ near $0$ is the graph of a complex valued function over some part of ${\bf C}\times\{0\}$. We conclude from (\ref{taylor-expansion-of-v}) that the set $\mbox{ pr}_1(v(\stackrel{\circ}{D^+}))$ either lies in the positive half--plane in ${\bf C}$ or in the negative half--plane depending on the sign of the number $c$. Note that we have ruled out the case $c<0$ by assumption. We write
\[
v(z)=(a(z),b(z))\ ,\ u_{\tau}(z)=(a_{\tau}(z),b_{\tau}(z))
\]
for suitable complex valued functions $a,a_{\tau},b,b_{\tau}$ defined on the upper half--disk $D^+$. These functions have the properties
\begin{itemize}
\item $b(0)=0$ , $Db(0)=0$ , $a(0)=0$ and $a(z)=c z+O(|z|^2)$ for some $c>0$,
\item $a_0(z)=z$ , $b_0\equiv 0$,
\item $b(D^+\cap{\bf R})\subset [0,\infty)\,,\,b_{\tau}(D^+\cap{\bf R})\subset (-\infty,0)$ if $\tau<0$,
\end{itemize}
(the last item in the list might as well be $b(D^+\cap{\bf R})\subset (-\infty,0]$, but then $b_{\tau}(D^+\cap {\bf R})\subset (0,+\infty)$ for $\tau<0$, which does not change the argument of the proof).
Restricting all functions involved to a smaller half--disk and confining ourselves to values of $\tau$ close to zero, we may assume that all the maps $a_{\tau}, a$ are local diffeomorphisms near $0$. Note that we did not parameterize $v$ so that it looks like $(z,\tilde{b}(z))$ since we want to keep the property $J(z,0)=i$. We compute now, denoting the partial derivative with respect to the second argument by $D_2J$,
\begin{eqnarray*}
0 & = & \pas v+J(v)\pat v \\
 & = & (\pas a,\pas b)+J(a,b)(\pat a,\pat b) \\
 & = & (\pas a,\pas b)+\left(J(a,0)+\int_0^1D_2J(a,\kappa b)\,b\,d\kappa\right)(\pat a,\pat b)\\
 & = & (\pas a+i\pat a+\alpha b,\pas b+i\pat b+\beta b),
\end{eqnarray*}
where we wrote
\[
(\alpha x,\beta x)=\left(\int_0^1D_2J(a,\kappa b)\,x\,d\kappa\right)(\pat a,\pat b).
\]
We consider the second component of the above differential equation
\[
0=\pas b+i\pat b+\beta b.
\]
The boundary version of the similarity principle applies here (see theorem \ref{similarity-principle-bdry}). If the $\infty$--jet of $b$ vanished at $0$ then we would have $b\equiv 0$ because of the similarity principle. This would imply that $v(z)=(a(z),0)$ (and $a$ is biholomorphic) and the assertion of the proposition follows.\\
Without the assumption $v(D^+_{\ve})\subset H^+\times{\bf C}$, the image of $v|_{D^+_{\ve}}$ would be contained in the image of the 'Schwarz reflection' of $u_0$.\\ 

We claim that the $\infty$--jet of $b$ actually has to vanish at $0$. Arguing indirectly, we assume that it does not. Without loss of generality we may also assume that $b(z)\neq 0$ for $z\in D^+\backslash\{0\}$ since $0$ is then an isolated intersection point of $v$ with $u_0$. Applying the similarity principle we find a holomorphic map $\sigma:D^+_{\varepsilon}\rightarrow{\bf C}$ defined on an upper half--disk of radius $0<\varepsilon<1$ and a map $\Phi\in\bigcap_{2<p<\infty}W^{1,p}(D^+_{\varepsilon},{\bf C}\backslash\{0\})$ with $\Phi((-\varepsilon,\varepsilon))\subset{\bf R}\backslash\{0\}$ so that
\[
b(z)=\Phi(z)\sigma(z)\ \mbox{on}\ D^+_{\varepsilon}.
\]
Write
\[
\sigma(z)=\sum^{\infty}_{k=k_0}a_kz^k
\]
with some $k_0\ge 2$. We know that $b|_{D^+\cap {\bf R}}$ does not change its sign near $0$ since $v$ does not intersect any of the solutions $u_{\tau}$ for $\tau<0$ (see figure \ref{boundary-curves}), therefore $k_0$ must be even. Recall that all the curves $u_{\tau}|_{D^+\cap {\bf R}}$ lie in the lower half--plane in ${\bf R}^2$ so that $b(D^+\cap {\bf R})\subset [0,\infty)$. We may assume, by making $\ve$ smaller if necessary, that all the maps $a_{\tau}$ are local diffeomorphisms if $|\tau|$ is sufficiently small since $a_0(z)=z$. We note that
\begin{equation}\label{no-solutions-on-the-boundary1}
-\delta\not\in b(\partial D^+_{\ve}),
\end{equation}
where $\partial D^+_{\ve}:=\{z\in D^+_{\ve}\,|\,z\in {\bf R}\ \mbox{or}\ |z|=\ve\}$, for all $\delta_0>\delta>0$ with some sufficiently small $\delta_0>0$. This is true because 
\[
0\not\in b(\overline D^+_{\ve}\backslash\{0\})
\]
and $b((-\ve,\ve))\subset [0,\infty)$ since the path $b|_{D^+\cap\{|z|=\ve\}}$ starts and ends somewhere on the positive real axis and avoids the origin. If $|\tau|$ is sufficiently small, we may also assume that
\begin{equation}\label{no-solutions-on-the-boundary2}	
-\delta\neq b(z)-b_{\tau}(a_{\tau}^{-1}(a(z)))\ ,\ \mbox{for all}\ z\in \partial D^+_{\ve}\ \mbox{and all}\ 0<\delta<\delta_0/2
\end{equation}
since $b_0\equiv 0$. The Brouwer degrees $\degree(b,D^+_{\ve},-\delta)$ and $\degree(b-b_{\tau}\circ a_{\tau}^{-1}\circ a,D^+_{\ve},-\delta)$ are then well--defined for all sufficiently small $\delta>0$ and they agree in view of (\ref{no-solutions-on-the-boundary1}) and (\ref{no-solutions-on-the-boundary2}). We observe that for all small $\delta$
\[
\degree(b,D^+_{\ve},-\delta)=\degree(\sigma,D^+_{\ve},-\delta)
\]
since $\Phi$ can be removed by a trivial homotopy argument. We continue now the holomorphic map $\sigma$ analytically onto the whole disk $D_{\ve}$ by Schwarz reflection. Choosing $\delta>0$ sufficiently small we may assume that $|\sigma(z)|>\delta$ for all $z\in\partial D_{\ve}$ so that the Brouwer degree $\degree(\sigma,D_{\ve},-t\delta)$ is defined for all $0\le t\le 1$. Then
\begin{eqnarray*}
k_0 & = & \degree(\sigma,D_{\ve},0) \\
 & = & \degree(\sigma, D_{\ve},-\delta) \\
 & = & \degree(\sigma, D^-_{\ve},-\delta)+\degree(\sigma,D^+_{\ve},-\delta)\\
 & = & 2\,\degree(\sigma,D^+_{\ve},-\delta)\\
 & = & 2\,\degree(b,D^+_{\ve},-\delta) \\
\end{eqnarray*}
and
\begin{equation}\label{degree-is-nontrivial}
\degree(b,D^+_{\ve},-\delta)=\frac{k_0}{2}\ge 1.
\end{equation}
By assumption the images of $v$ and $u_{\tau}$ are disjoint if $\tau<0$, i.e. the equations
\[
a(z)=a_{\tau}(z)\ ,\ b(z)=b_{\tau}(z),
\]
have no solution if $\tau<0$. We rewrite them as follows:
\[
b(z)=b_{\tau}(a_{\tau}^{-1}(a(z)))\ ;\ z\in D^+_{\ve}.
\]
Since there is no solution to this equation we have necessarily 
\begin{eqnarray*}
0 & = & \degree(b-b_{\tau}\circ a_{\tau}^{-1}\circ a,D^+_{\ve},0)\\
 & = & \degree(b-b_{\tau}\circ a_{\tau}^{-1}\circ a,D^+_{\ve},-t\delta)\,,0\le t\le 1\\
 & = & \degree(b-b_{\tau}\circ a_{\tau}^{-1}\circ a,D^+_{\ve},-\delta)
\end{eqnarray*}
in contradiction to (\ref{degree-is-nontrivial}). This completes the proof. 

\qed

\begin{figure}[!ht]
\begin{center}
\epsfig{file=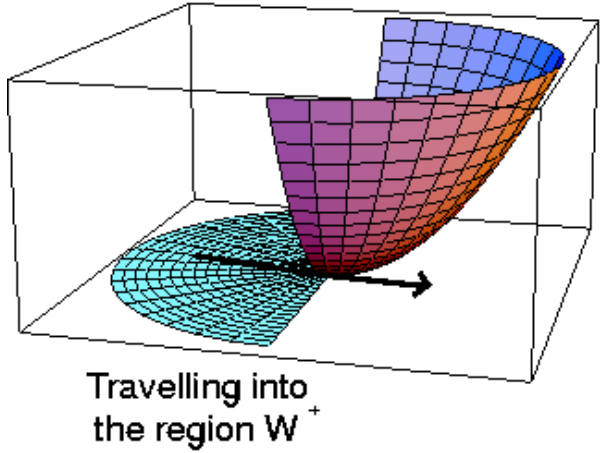,height=6cm,angle=0}
\caption{Following the arrow from the graph of $u_0$ onto the graph of $v$ would lead into the region $W^+$.}\label{bad-cannot-happen}
\end{center}
\end{figure}

{\bf Proof of theorem \ref{interior-and-boundary-intersections}:}\\
Assume first that $\tu_0(p)=(a_0(p),u_0(p))=\tilde{v}(q)$ for some $p\in{\bf R}\times (0,1]$. Let us first consider the case $p\in{\bf R}\times\{1\}$ so that $a_0(p)=0$. If $\tilde{v}$ is also a pseudoholomorphic strip then proposition \ref{maximum-principle} implies that $q\in{\bf R}\times \{1\}$ as well. If $\tilde{v}=(b,v)$ is a pseudoholomorphic disk then $b$ has to be negative on the interior of the disk since $\Delta b\ge 0$ with zero boundary conditions. Hence $q\in\partial D$. If $\tilde{v}$ is a pseudoholomorphic cylinder then trivially $q\in\partial Z^-$. After invoking lemma \ref{local-coordinates-near-boundary-point} we may view both $\tu_0$ and $\tilde{v}$ as maps on the half--disk $D^+$ into ${\bf C}^2$ with boundary conditions $\tu_0(D^+\cap {\bf R})\,,\,\tilde{v}(D^+\cap {\bf R})\subset {\bf R}^2$. The proof of the case $p\in{\bf R}\times\{1\}$ would be complete if we could show that $\mbox{pr}_1[\tilde{v}(D^+_{\ve})]$ lies in the upper--half plane of ${\bf C}$ (see figure \ref{good}) with $\mbox{pr}_1:{\bf C}^2\rightarrow{\bf C}$ being the projection onto the first factor. We may write ${\bf R}\times M$ as the disjoint union of the two domains ${\bf R}^{\pm}\times M$ and the hypersurface $\{0\}\times M$. We note that $\tu_0$ and $\tilde{v}$ map the interiors of their domains into ${\bf R}^-\times M$, and the boundaries ${\bf R}\times\{1\}$, $\partial D$ and $\partial Z^-$ respectively into the hypersurface $\{0\}\times M$ (which is called pseudoconvex with respect to ${\bf R}^-\times M$). In the local picture in ${\bf C}^2$ we obtain two domains $W^{\pm}$ corresponding to ${\bf R}^{\pm}\times M$ and a real hypersurface containing ${\bf R}^2$. If $\gamma(t)$ is any differentiable path in ${\bf C}^2$ which agrees with $\tu_0(0,1-t)$ for $0\le t\le 1$ then $\gamma(1+\ve)\in W^+$ for all small $\ve>0$. This holds because of proposition \ref{maximum-principle} since the outer normal derivative of $a_0$ in the point $(0,0)$ is strictly positive. 
If $\mbox{pr}_1[\tilde{v}(D^+_{\ve})]$ were contained in the lower half--plane in ${\bf C}$ then we could find a differentiable path $\gamma(t)$ as above which lies in the image of $\tilde{v}$ for $t\ge 1$, $t-1$ small (see figure \ref{bad-cannot-happen}), because $\tu_0$ and $\tilde{v}$ are tangent at $(0,0)$. This would imply that $\tilde{v}$ maps some interior points of its domain into $W^+$ which is impossible.
Hence we have reduced the case $p\in{\bf R}\times \{1\}$ to the local proposition \ref{intersection-of-families-local}. Assume now that $p\in{\bf R}\times (0,1)$, i.e. $a_0(p)<0$ by proposition \ref{maximum-principle}. Then the point $q$ is contained in ${\bf R}\times[0,1)$, $\stackrel{\circ}{D}$ or $(-\infty,0)\times S^1$. If $q$ is in ${\bf R}\times(0,1)$, $\stackrel{\circ}{D}$ or $(-\infty,0)\times S^1$ then we can apply proposition \ref{intersection-of-families-local-interior} right away. If $q\in{\bf R}\times\{0\}$ then we can reflect $\tilde{v}$ locally near $q$ at the bondary using lemma \ref{local-coordinates-near-boundary-point} and apply proposition \ref{intersection-of-families-local-interior}. The issue of $\tilde{v}(D^+_{\ve})$ lying on the wrong side of ${\bf R}^2$ does not come up here because $p$ is an interior point.\\  
We are left with the case where $p\in{\bf R}\times\{0\}$. The boundaries of Bishop disks and periodic orbits of the Reeb vector field cannot intersect ${\mathcal L}=\partial {\mathcal D}$, so we only deal with the case where $\tilde{v}$ is another pseudoholomorphic strip. The case $q\in{\bf R}\times\{1\}$ can not occur because of proposition \ref{maximum-principle}. If $q$ is an interior point then we may reflect $\tu_0$ locally near $p$ and apply proposition \ref{intersection-of-families-local-interior}. If now both $q,p$ are boundary points on ${\bf R}\times\{0\}$ then unfortunately, the boundary condition ${\bf R}\times {\mathcal L}$ is not contained in a pseudoconvex hypersurface in ${\bf R}\times M$, as it was the case with $\{0\}\times {\mathcal D}^{\ast}$. Hence the image of $\tilde{v}$ may either agree with the image of $\tu_0$ or with image of the reflection of $\tu_0$. Writing $v=(a,\theta,x,y)$ and $u_0=(a_0,\theta_0,x_0,y_0)$ in local coordinates near $\{0\}\times{\mathcal L}$ we use lemma \ref{local-coordinates-near-boundary-point} in conjunction with proposition \ref{reflection-2}. The first factor in ${\bf C}^2$ corresponds to the $\theta x$--plane. The condition whether $\mbox{pr}_1[\tilde{v}(D^+_{\ve})]$ lies in the upper or the lower half--plane is the same as asking whether $\pat x(q)=-\pas\theta(q)$ and $\pat x_0(p)=-\pas\theta_0(p)$ have the same sign or not. Hence the question is whether $\pas u_0(p),\pas v(q)\in T_{v(q)}{\mathcal L}$ have the same orientation or not.

\qed

We study now the situation where the first intersection occurs at infinity.

\begin{theorem}\label{intersections-at-infinity}
Let $\tu_{\tau}=(a_{\tau},u_{\tau})$ , $-1<\tau\le 0$ be a smooth family of embedded solutions of (\ref{main-boundary-value-problem}) with pairwise disjoint images. Let $\tilde{v}$ be another embedded solution. Assume that all the maps $\tu_{\tau}$ and $\tilde{v}$ have the same exponential decay rate $\lambda_+=-\pi/2$ as $s\rightarrow+\infty$. We assume also that $\tilde{v}$ and $\tu_0$ converge to the same point on $\{0\}\times{\mathcal L}$ as $s\rightarrow+\infty$, but 
\[
\tilde{u}_{\tau}(S)\cap \tilde{v}(S)=\emptyset\ \mbox{for all}\ \tau<0
\]
Then 
\[
\tilde{v}({\bf R}\times\{0\})=\tu_{0}({\bf R}\times\{0\}).
\] 
\end{theorem}
{\bf Proof:}\\
The idea is to reduce the situation to the one in proposition \ref{intersection-of-families-local-interior}. Since we are only interested in large $s$ we may work in local coordinates near $\{0\}\times{\mathcal L}$. Using proposition \ref{reflection-2} we reflect all the solutions at the boundary $[s_0,\infty)\times\{0\}$, and we obtain pseudoholomorphic strips defined on $S_+=[s_0,\infty)\times[-1,+1]$ having boundary condition on $\{0\}\times{\mathcal D}$. We consider now the biholomorphic map
\[
\phi:S_+\longrightarrow \Omega\backslash\{0\}
\]
\[
\phi(z)=w=e^{-\frac{\pi}{2}z},
\]
where $\Omega=\{z\in{\bf C}\,|\,\mbox{Re}(z)\ge 0\}$.
Composing all the reflected pseudoholomorphic strips $\tilde{v}\,,\tu_{\tau}$ with 
\[
\psi(w)=\phi^{-1}(w)=-\frac{2}{\pi}\log\,w
\]
we obtain pseudoholomorphic curves on the punctured half--disk with boundary condition in $\{0\}\times{\mathcal D}$. Our aim is to reflect them once again at the boundary near the origin using corollary \ref{local-reflection}. It is clear that all curves extend continuously over the origin. We have to verify that their derivatives also extend and that they are not zero in the origin. This is where the decay rate of $-\frac{\pi}{2}$ comes in. If we carried out the reflection and reparametrization procedure with a pseudoholomorphic strip of faster decay, we would obtain a pseudoholomorphic half--disk with vanishing derivative in the origin. Using the asymptotic formula (\ref{odd-explicit-asymptotic-formula-v2}) and identifying ${\bf R}^4$ with ${\bf C}^2$ we may write the resulting pseudoholomorphic half--disks as follows:
\[
\tu_0(\psi(w))=\rho(w)\frac{w}{|w|}\,(1,q(0))+\rho(w)\varepsilon(\psi(w)),
\]
where
\[
\rho(w)=\kappa_0\,\exp\left(\int_{s_0}^{-\frac{2}{\pi}\log|w|}\alpha(\tau)\,d\tau\right)\ ,\ \kappa_0>0
\]
and $c>0$ is some constant. We note that $D^{\alpha}\varepsilon(\psi(0))=0$ for all $|\alpha|\ge 0$. We compute with $\partial=\partial_w$, $\bar{\partial}=\partial_{\bar{w}}$:
\[
\bar{\partial}\rho(w)=-\frac{w}{\pi|w|^2}\rho(w)\ \alpha\left(-\frac{2}{\pi}\log|w|\right),
\]
\[
\partial\rho(w)=-\frac{\bar{w}}{\pi|w|^2}\rho(w)\ \alpha\left(-\frac{2}{\pi}\log|w|\right),
\]
which are both bounded since $|\rho(w)|\le\,c\,|w|$ for some constant $c>0$. We note also that
\[
\rho(w)=\kappa_0\,e^{\frac{\pi}{2}s_0}\,|w|\,\exp\left(\int_{s_0}^{-\frac{2}{\pi}\log|w|}[\alpha(\tau)+\frac{\pi}{2}]d\tau\right).
\]
Then
\begin{eqnarray*}
\partial\left(\rho(w)\frac{w}{|w|}\right) & = & \frac{\rho(w)}{2|w|}-\frac{\bar{w}}{\pi|w|^2}\rho(w)\alpha\left(-\frac{2}{\pi}\log|w|\right)\cdot\frac{w}{|w|}\\
 & = & \frac{\rho(s)}{|w|}-\frac{\rho(s)}{\pi|w|}\left(\frac{\pi}{2}+\alpha\left(-\frac{2}{\pi}\log|w|\right)\right)\\
 & = & \kappa_0\,e^{\frac{\pi}{2}s_0}\,\exp\left(\int_{s_0}^{-\frac{2}{\pi}\log|w|}[\alpha(\tau)+\frac{\pi}{2}]d\tau\right)-\\
  & & -\frac{\rho(s)}{\pi|w|}\left(\frac{\pi}{2}+\alpha\left(-\frac{2}{\pi}\log|w|\right)\right).
\end{eqnarray*}
The first term equals
\[
\kappa_0\,e^{\frac{\pi}{2}s_0}\,\exp\left(\int_{s_0}^{+\infty}[\alpha(\tau)+\frac{\pi}{2}]d\tau\right)\neq 0
\]
for $w=0$ by theorem \ref{convergence-of-alpha}, and the second one converges to zero since $|\rho(w)|\le c|w|$ and $\alpha(s)\rightarrow -\frac{\pi}{2}$ as $s\rightarrow+\infty$. We evaluate
\[
\bar{\partial}\left(\rho(w)\frac{w}{|w|}\right)=-\rho(w)\frac{w^2}{\pi|w|^3}\left(\frac{\pi}{2}+\alpha\left(-\frac{2}{\pi}\log|w|\right)\right),
\]
which converges to zero as $w\rightarrow 0$. We compute
\begin{eqnarray*}
|\rho(w)\bar{\partial}(\ve\circ\psi)(w)| & \le &  c|w|\,|D\ve(\psi(w))|\,|\partial\psi(w)|\\
 & \le & c\,|w|\,|D\ve(\psi(w))|\,| w|^{-1}\stackrel{|w|\rightarrow 0}{\longrightarrow} 0.
\end{eqnarray*}
We estimate with $\delta\ge\frac{\pi}{2}$ as in theorem \ref{convergence-of-alpha} 
\[
|\ve(\psi(w))|\le\,c\,e^{\frac{2\delta}{\pi}\log|w|}=c\,|w|^{\frac{2\delta}{\pi}}
\]
and 
\[
|\bar{\partial}\rho(w)\,(\varepsilon\circ\psi)(w)|\le c\,|w|^{\frac{2\delta}{\pi}}
\]
which also tends to zero as $|w|\rightarrow 0$. We proceed similarly for $\partial(\rho(\varepsilon\circ\psi))$ which also vanishes for $w=0$. Hence we may view all the pseudoholomorphic strips $\tu_{\tau},\tilde{v}$ for large $s$ as embedded pseudoholomorphic disks centered at the origin because their derivatives are not zero in the origin. Moreover, $\tilde{v}$ and $\tu_0$ have an isolated intersection in the origin, while $\tilde{v}$ does not intersect any of the disks $\tu_{\tau}$ for $\tau<0$. This completes the proof of the theorem. \qed
 
\section{Remarks about the implicit function theorem and about transversality}\label{IFT-second-version}
The main result of the paper \cite{part2} is theorem \ref{main-implicit-function-theorem}, the implicit function theorem. We have assumed that the solution $\tu_0$ decays to its endpoints at exponential rates $\exp(-|\lambda_{\pm}s|)$ with $|\lambda_{\pm}|\le\pi$, and that the Maslov--index ot $\tu_0$ vanishes. The purpose of this section is to show that there is also some version of theorem \ref{main-implicit-function-theorem} for arbitrary decay rates at the ends if the Maslov index assumes suitable values. In this section we are using the notation of the paper \cite{part2}, and we will indicate the necessary modifications in order to establish such an implicit function theorem.\\
In the paper \cite{part2} we started with an embedded solution $\tu_0$ to the boundary value problem \ref{main-boundary-value-problem}, and we attempted to find more solutions $\tu$ near by of the form 
\begin{equation}\label{normal-solution}
\tu(s,t)=\Phi_{c_-,c_+}(s,t,x(s,t),y(s,t)),
\end{equation}
where
\[
\Phi_{c_-,c_+}:S\times{\bf R}^2\supset S\times B_{\ve}(0)\longrightarrow {\bf R}\times M\ ,\ c_-,c_+\in {\bf R}
\]
\begin{eqnarray*}
\Phi_{c_-,c_+}(s,t,x,y) & := & \exp_{\tu_0(s,t)}\Big(x\,n(s,t)+y\,m(s,t)+\\
 & & +[c_-\beta_-(s)+c_+\beta_+(s)]\,(0,1,0,0)\Big),
\end{eqnarray*} 
and where $n(s,t)\in T_{\tu_0(s,t)}({\bf R}\times M)$ and $m(s,t)=\tilde{J}(\tu_0(s,t))n(s,t)$ are suitable normal vectors to the solution $\tu_0$. The letters $\beta_{\pm}$ stand for cut--off functions which equal $1$ for $|s|$ large and positive (or negative in the case of $\beta_-$). Their purpose is to move the end points of the map $\tu$ away from the end points of the solution $\tu_0$. In this section we will only require $\tu_0$ to be immersed. The implicit function theorem from \cite{part2} also applies to $\tu_0$ which are merely immersed, but the near by solution will then of course also be immersed only and not necessarily embedded. In the paper \cite{part2} we have set up a partial differential equation for the map $(s,t)\mapsto(x(s,t),y(s,t))\in {\bf R}^2$ so that $\tu$ given by (\ref{normal-solution}) is a $\tilde{J}$--holomorphic curve. The investigation of this PDE uses the following weighted Sobolev spaces 
\[
H^{2,p,\gamma}_L(S,{\bf C}):=\{u\in H^{2,p}(S,{\bf C})\,|\,\|u\|_{2,p,\gamma}:=\|\rho^{\gamma}u\|_{2,p}<\infty\ ,
\]
\[
u(s,0)\in{\bf R}\ ,\ u(s,1)\in{\bf R}\cdot (a_1(s)+ia_2(s))\ \},
\]
\[
H^{1,p,\gamma}(S,{\bf C}):=\{u\in H^{1,p}(S,{\bf C})\,|\,\|u\|_{1,p,\gamma}:=\|\rho^{\gamma}u\|_{1,p}<\infty\ \},
\]
where $\gamma:{\bf R}\rightarrow{\bf R}$ is a smooth function with
\[
\gamma(s)\stackrel{s\rightarrow\pm\infty}{\longrightarrow}\gamma_{\pm}\ ,\ \frac{d^k}{ds^k}\gamma(s)\stackrel{s\rightarrow\pm\infty}{\longrightarrow}0\ ,\ k\ge 1,
\]
and where $\rho(s)$ is a smooth function which agrees with $e^{\int_{s_0}^s\alpha_{\pm}(\tau)d\tau}$ for large $|s|$, and $\alpha_{\pm}$ are the functions appearing in the asymptotic formula for $\tu_0$ (theorem \ref{asymptotic-formula-theorem}). Recall also that $\alpha_{\pm}(s)\rightarrow\lambda_{\pm}$ as $s\rightarrow\pm\infty$. The map $s\mapsto a_1(s)+i a_2(s)\in {\bf C}\backslash \{0\}$ takes care of the boundary condition (see definition \ref{admissible-normal-vector} below for details). Since $p>2$, the above spaces consist of differentiable and continuous functions respectively. If $\gamma_{\pm}<0$ then the above Sobolev spaces consists of functions with a certain exponential decay at infinity. The PDE for $(x,y)$ makes sense under the following assumptions on the weights:
\begin{itemize}
\item If $\lambda_{\pm}=\mp\frac{\pi}{2}$ and $-\frac{1}{2}<\gamma_{\pm}<0$, or
\item If $\lambda_{\pm}\in{\bf Z}\pi$ and $-\frac{1}{2}-\frac{\delta}{|\lambda_{\pm}|}<\gamma_{\pm}<-\frac{1}{2}$, where $\delta>0$ is the exponential decay rate of the remainder terms in the asymptotic formula (theorem \ref{convergence-of-alpha}).
\end{itemize}
We will investigate the following question: Given an immersed solution $\tu_0$ with decay rates $\lambda_{\pm}$ near the ends and with Maslov index $\mu(\tu_0)$, when are there solutions $\tu$ near by of the form (\ref{normal-solution}) with $c_-=c_+=0$ ? This means that $\tu$ will have the same end points as the solution $\tu_0$ we started with. Let us recall the definition of the Maslov index. 
\begin{definition}\label{admissible-normal-vector}
If $\tu$ is an immersed solution to the boundary value problem (\ref{main-boundary-value-problem}) then we call a section $n$ in $\tu^{\ast}T({\bf R}\times M)$ an admissible normal vector if it satisfies the following conditions:
\begin{itemize}
\item $n(s,t)\not\in\,\mbox{Span}\{\pas\tu(s,t),\pat\tu(s,t)\}$,
\item $n(s,0)\in T_{\tu(s,0)}({\bf R}\times{\mathcal L})$,
\item $n(s,t)\longrightarrow n_{\pm\infty}(t)$ uniformly in $t$ as $s\rightarrow \pm\infty$, where
\[
n_{\pm\infty}(t):=\left\{\begin{array}{ccc} (0,\pm 1,0,0) & \mbox{if} & \lambda_{\pm}=\mp\frac{\pi}{2}\\ (\cos\frac{\pi t}{2},-q_{\pm}(0)\cos\frac{\pi t}{2},0,\mp\sin\frac{\pi t}{2}) & \mbox{if} & \lambda_{\pm}\in{\bf Z}\pi\end{array}\right.
\]
and $\lambda_{\pm}$ as in theorem \ref{asymptotic-formula-theorem}.
\item There is a path 
\[
\gamma=a_1+ia_2:{\bf R}\rightarrow {\bf C}\backslash\{0\}
\]
so that 
\[
a_1(s)\,n(s,1)+a_2(s)\,\tilde{J}(\tu(s,1))n(s,1)\in T_{\tu(s,1)}(\{0\}\times{\mathcal D}).
\]
\end{itemize}
Writing $\gamma(s)=\rho(s)\,e^{i\phi(s)}$ and $\phi_{\pm}=\lim_{s\rightarrow\pm\infty}\phi(s)$, the Maslov index of $\tu$ is then defined by 
\[
\mu(\tu):=\frac{1}{\pi}(\phi_+-\phi_-).
\]
\end{definition}
We note that $n_{\pm\infty}(t)$ satisfies the boundary conditions $n_{\pm\infty}(0)\in T_{\tu(s,0)}({\bf R}\times{\mathcal L})$ and $n_{\pm\infty}(1)\in T_{\tu(s,1)}(\{0\}\times{\mathcal D})$ for $|s|$ large. Theorem \ref{main-implicit-function-theorem} assumes that the Maslov--index of $\tu_{0}$ is zero and that $\tu_0$ has particular rates of decay, but the computation of the Fredholm index in the paper \cite{part2} was actually carried out without these assumptions. Proposition 4.1 in \cite{part2} establishes the Fredholm property and the index computation based on computing the spectral flow of a Cauchy Riemann type operator
\[
T:H^{1,p}_{{\bf R}}(S,{\bf C})\longrightarrow L^p(S,{\bf C})
\]
\[
(T\eta)(s,t):=\pas\eta(s,t)+i\pat\eta(s,t)+F(s)\eta(s,t),
\]
where $F$ is a monotonous smooth function with
\[
F(s)\rightarrow\left\{\begin{array}{ccc} -\gamma_-\lambda_- & \mbox{as} & s\rightarrow-\infty\\
-\gamma_+\lambda_+-\pi\mu(\tu_{\tau_0}) & \mbox{as} & s\rightarrow+\infty\end{array}\right.
\]
where $\lambda_{\pm}$ are the decay rates of $\tu_{\tau_0}$ in the asymptotic formula and
\[
H^{1,p}_{{\bf R}}(S,{\bf C}):=\{\eta\in H^{1,p}_{{\bf R}}(S,{\bf C})\,|\,\eta(\partial S)\subset{\bf R}\}.
\]
The proof of proposition 4.1 in \cite{part2} does not use the assumptions $\lambda_{\pm}\in\{\mp\pi,\mp\frac{\pi}{2}\}$ and $\mu(\tu_{0})=0$ until the very last line of the proof where a formula is given for the Fredholm index. Let $A(s)$ be the operator $\gamma\mapsto -i\dot{\gamma}-F\gamma$ acting on $H^{1,2}_{{\bf R}}([0,1],{\bf C})$, and denote their eigenvalues by $\lambda_n(s)=n\pi-F(s)$. Then
\begin{eqnarray*}
 & & \mbox{ind}(T)\\ & = & -\sum_{\{(n_0,s_0):\lambda_{n_0}(s_0)=0\}}\mbox{sign}\,\lambda'_{n_0}(s_0)\\
 & = & \left\{\begin{array}{cc} \#({\bf Z}\pi\cap[\gamma_+\lambda_++\pi\mu(\tu_0),\gamma_-\lambda_-]) & \mbox{if}\ \gamma_+\lambda_++\pi\mu(\tu_0)<\gamma_-\lambda_-\\
  -\#({\bf Z}\pi\cap [\gamma_-\lambda_-,\gamma_+\lambda_++\pi\mu(\tu_0)]) & \mbox{if}\ \gamma_-\lambda_-<\gamma_+\lambda_++\pi\mu(\tu_0)
  \end{array}\right.
\end{eqnarray*}
If $\lambda_{\pm}=\mp m_{\pm}\pi$ for some positive integers $m_{\pm}$ and if we have chosen the weights such that $\gamma_-=\gamma_+=\gamma$ then $\mbox{ind}(T)$ is positive if and only if $\mu(\tu_0)<-|\gamma|\,(m_-+m_+)$. A sufficient condition for $\mbox{ind}(T)>0$ is 
\begin{equation}\label{index-is-positive-1}
\mu(\tu_0)<-\frac{1}{2}(m_-+m_+).
\end{equation}
We write $m_++m_-=2M$ or $m_-+m_+=2M+1$ for a suitable integer $M\ge 1$. We then obtain under the assumption (\ref{index-is-positive-1})
\begin{eqnarray*}
\mbox{ind}(T) & = & \#({\bf Z}\cap[-\gamma m_++\mu(\tu_0),\gamma m_-])\\
 & = & \left\{\begin{array}{cc} \#({\bf Z}\cap [-2M\gamma+\mu(\tu_0),0]) & \mbox{if}\ m_++m_-=2M\\
 \#({\bf Z}\cap[-\gamma(2M+1)+\mu(\tu_0),0]) & \mbox{if}\ m_++m_-=2M+1\end{array}\right.\\
  & = & -(M+\mu(\tu_0)).
\end{eqnarray*} 
If the decay rate at one end is $\mp\frac{\pi}{2}$ and at the other an integer multiple $\pm m\pi$ of $\pi$ with $m\ge 1$ then a sufficient condition for $\mbox{ind}(T)$ being positive is 
\begin{equation}\label{index-is-positive-2}
\mu(\tu_0)<-\frac{1}{2}m-\frac{1}{4}.
\end{equation}
Assuming that the condition (\ref{index-is-positive-2}) is satisfied we compute
\begin{eqnarray*}
\mbox{ind}(T) & = & \left\{\begin{array}{cc} -(\mu(\tu_0)+\frac{m}{2}) & \mbox{if $m$ is even}\\
-(\mu(\tu_0)+\frac{m-1}{2}) & \mbox{if $m$ is odd}\end{array}\right.
\end{eqnarray*}
The transversality argument is based on investigating the kernel of an operator similar to $T$ above with the same spectral flow (we denote it again by $T$). The proofs in the paper \cite{part2} also work in the more general context after some minor changes which we will indicate now. References refer to formulae and results in the paper \cite{part2}. Formula (4.12) in \cite{part2} is valid too with the modification
\[
\Gamma(s,t)\rightarrow\left\{\begin{array}{ccc} -\gamma_-\lambda_- & \mbox{as} & s\rightarrow-\infty\\ -\gamma_+\lambda_+-\pi\mu(\tu_{\tau_0}) & \mbox{as} & s\rightarrow+\infty\end{array}\right.
\]
In lemma 4.4 we should replace formula (4.16) by
\[
A_+=-i\frac{d}{dt}+(\gamma_+\lambda_++\pi\mu(\tu_{\tau_0}))\,\mbox{Id}
\]
and the eigenvalues of the operators $A_{\pm}$ are given by
\[
\nu_n^+=n\pi+\gamma_+\lambda_++\pi\mu(\tu_{\tau_0})
\]
and
\[
\nu_n^-=n\pi+\gamma_-\lambda_-.
\]
while the formula (4.17) for the corresponding eigenvectors remains unchanged
\[
e_{\pm}(t)=e^{i\,n\pi t}.
\]
The assertion of proposition 4.6 in \cite{part2} was the following formula relating the numbers of the zeros of any nontrivial element $\eta$ in the kernel of $T$ to the asymptotic eigenvalues $\nu_{\pm}^n$
\begin{equation}\label{formel-fuer-nullstellen}
n_+-n_-=2\sum_{z\in N_{int}}o(z)+\sum_{z\in N_{bd}}o(z)\ge 0.
\end{equation}
Here $N_{int}=\{z\in\stackrel{\circ}{S}\,|\,\eta(z)=0\}$ and $N_{bd}=\{z\in\partial S\,|\,\eta(z)=0\}$.
The inequality $n_+-n_-<0$ which was derived from $\nu_n^+<0$ and $\nu_n^->0$ in the paper \cite{part2} has to be replaced by the inequalities
\[
\nu_+=\pi n_++\gamma_+\lambda_++\pi\mu(\tu_{\tau_0})<0\ ,\ \nu_-=\pi n_-+\gamma_-\lambda_->0
\]
which leads to
\begin{equation}\label{trans-inequality}
n_+-n_-<\frac{1}{\pi}(\gamma_-\lambda_--\pi\mu(\tu_0)-\gamma_+\lambda_+).
\end{equation}
In the case where $\lambda_{\pm}=\mp m_{\pm}\pi$ and $\gamma_-=\gamma_+=\gamma$ with $\gamma<-\frac{1}{2}$ we obtain
\begin{eqnarray*}
n_+-n_- & < & \gamma m_--\mu(\tu_0)+\gamma m_+\\
 & < & \left\{\begin{array}{cc} -M-\mu(\tu_0) & \mbox{if}\ m_++m_-=2M\\
 -M-\mu(\tu_0)-\frac{1}{2} & \mbox{if}\ m_++m_-=2M+1\end{array}\right..
\end{eqnarray*}
This implies
\begin{equation}\label{trans-2}
n_+-n_-\le \mbox{ind}(T)-1
\end{equation}
because $n_+-n_-$ is an integer. Consider now the case where the decay rate in one end, say the positive end, is $\lambda_+=-\frac{\pi}{2}$ and the decay rate in the other is $\lambda_-=m\pi$, $m\in{\bf N}$. The weights satisfy the conditions $\gamma_+<0$ and $\gamma_-<-\frac{1}{2}$. Then 
\begin{eqnarray*}
n_+-n_- & < & \gamma_-m-\mu(\tu_0)+\frac{\gamma_+}{2}\\
 & < & -\frac{m}{2}-\mu(\tu_0)\\
 & = & \left\{\begin{array}{cc} \mbox{ind}(T) & \mbox{if $m$ is even}\\
 \mbox{ind}(T)-\frac{1}{2} & \mbox{if $m$ is odd}\end{array}\right.
\end{eqnarray*}
and we conclude again
\[
n_+-n_-\le\mbox{ind}(T)-1.
\]
In view of equation (\ref{formel-fuer-nullstellen}) we conclude that any nontrivial element $\eta\in\ker T$ has at most $(\mbox{ind}(T)-1)$ zeros on $\partial S$. We claim that this implies that the operator $T$ must be surjective. Arguing indirectly, we assume that $T$ has a nontrivial cokernel. Then we must have 
\[
\Lambda=\mbox{dim}\,\ker T\ge\,\mbox{ind}(T)+1.
\]
We pick now linear independent elements $\eta_1,\ldots,\eta_{\Lambda}\in\ker T$ and points $z_1,\ldots,z_{\mbox{ind}(T)}\in\partial S$ such that 
\[
{\bf R}\ni\eta_k(z_l)\neq 0\ \forall\ 1\le k\le\Lambda\ ,\ 1\le l\le \mbox{ind}(T).
\]
The expressions
\[
\sum_{k=1}^{\Lambda}\lambda_k\,\eta_k(z_l)=0\ ,\ 1\le l\le \mbox{ind}(T)
\]
make up a system of $\mbox{ind}(T)$ linear equations in $\Lambda$ variables with real coefficients since $\eta(\partial S)\subset {\bf R}$. Hence there is a nontrivial solution $(\lambda_1,\ldots,\lambda_{\Lambda})$. Because the $\eta_k$ were linear independent, we have constructed a nontrivial element $\sum_{k=1}^{\Lambda}\lambda_k\eta_k\in\ker T$ which has $\mbox{ind}(T)$ zeros on $\partial S$, a contradiction to equations (\ref{formel-fuer-nullstellen}) and (\ref{trans-2}). We summarize our discussion as follows:

\begin{theorem}\label{extended-implicit-function-theorem}
{\bf (Implicit Function Theorem--second version)}\\
Let $\tu_0=(a_0,u_0)$ be an immersed solution of (\ref{main-boundary-value-problem}). Denote by $\lambda_{\pm}$ the decay rates of $\tu_0$ as in theorem \ref{asymptotic-formula-theorem} and let $\mu(\tu_0)$ be the Maslov index of $\tu_0$. Assume that one of the following conditions are satisfied:
\begin{enumerate}
\item $\lambda_{\pm}=\mp m_{\pm}\pi$ with integers $m_{\pm}\ge 1$ and $\mu(\tu_0)<-\frac{1}{2}(m_-+m_+)$ holds,
\item the absolute value of one of the numbers $\lambda_{\pm}$ equals $\frac{\pi}{2}$, and the absolute value of the other equals $m\pi$ with some positive integer $m$. We also assume that $\mu(\tu_0)<-\frac{1}{2}m-\frac{1}{4}$.
\end{enumerate}
Assume moreover that $\mbox{dist}(u_0({\bf R}\times\{1\}),\Gamma)>0$. Then there is an integer $N\ge 1$ and a smooth family $(\tilde{v}_{\tau})_{\tau\in {\bf R}^N}$ of solutions the boundary value problem (\ref{main-boundary-value-problem}) with the following properties:
\begin{itemize}
\item $\tilde{v}_0=\tu_0$ and $\tilde{v}_{\tau}\not\equiv\tu_0$ if $\tau\neq 0$
\item the solutions $\tilde{v}_{\tau}$ have the same end points as the solution $\tu_0$, i.e.
\[
\lim_{s\rightarrow\pm\infty}\tu_0(s,t)=\lim_{s\rightarrow\pm\infty}\tilde{v}_{\tau}(s,t)\ \forall\ \tau,
\]
\item The solutions $\tilde{v}_{\tau}$ have the same Maslov--index and the same decay rates as $\tu_0$,
\end{itemize}
\end{theorem}
We omit the proof of the last assertion of theorem \ref{extended-implicit-function-theorem} since it is similar to the corresponding statement in \cite{part2}. 

\section{Proof of theorem \ref{compactness-result}}
Let $\tau_k$ be a sequence converging to $\tau_0$. The uniform gradient bound enables us to use regularity estimates and the Arzela--Ascoli theorem. They guarantee the existence of a subsequence, which we will denote again by $\tau_k$, such that the sequence $\tu_{\tau_k}$ converges in $C^{\infty}_{loc}$ to some map $\tu_{\tau_0}:S\rightarrow \RM$ which satisfies the differential equation $\pas\tu_{\tau_0}+\tilde{J}(\tu_{\tau_0})\pat\tu_{\tau_0}=0$, the boundary conditions $\tu_{\tau_0}(s,0)\in{\bf R}\times{\mathcal L}$ , $\tu_{\tau_0}(s,1)\in\{0\}\times{\mathcal D}^{\ast}$ and the condition
\[
\mbox{dist}\big(\{u_{\tau_0}(s,1)\,|\,s\in{\bf R}\}\,,\,\Gamma\big)>0.
\]
Because $u_{\tau}(0,0)=e$ for all $0\le\tau<\tau_0$ we also have $u_{\tau_0}(0,0)=e$. This ensures that $\tilde{u}_{\tau_0}$ is not constant (recall that the convergence is only in $C^{\infty}$ on compact sets; strips with two different ends might well converge in $C^{\infty}_{loc}$ to a constant map).
Let us show that the energy of $\tu_{\tau_0}$ is finite. We have shown in the paper \cite{part2} (proposition 2.3) that the energy of all the maps $\tu_{\tau_k}$ is bounded by some constant $\mbox{vol}_{\lambda}({\mathcal D})$ depending only on the Seifert surface ${\mathcal D}$ and the contact form $\lambda$. The condition in \cite{part2} that the path $s\mapsto\tu_{\tau_k}(s,1)$ on the Seifert surface represents a trivial homology class in $H_1({\mathcal D},\partial{\mathcal D})$ is of course satisfied here. In particular, for any compact subset $K\subset S$ and any $k$
\[
\sup_{\phi\in\Sigma}\int_K\tu_{\tau_k}^{\ast}\,d(\phi\lambda)\le \,\mbox{vol}_{\lambda}({\mathcal D}).
\]
We may pass to the limit $k\rightarrow\infty$ and obtain
\[
\sup_{\phi\in\Sigma}\int_K\tu_{\tau_0}^{\ast}\,d(\phi\lambda)\le \,\mbox{vol}_{\lambda}({\mathcal D})
\]
which implies $E(\tu_{\tau_0})\le \,\mbox{vol}_{\lambda}({\mathcal D})$. Finiteness of energy has many consequences, the most important one is that all the results about asymptotic behavior apply now to $\tu_{\tau_0}$. In particular, $\tu_{\tau_0}(s,t)$ converges to points $\tilde{p}_{\pm}$ on $\{0\}\times{\mathcal L}$ as $|s|\rightarrow\infty$ uniformly in $t$. Since the convergence $\tu_{\tau_k}\rightarrow\tu_{\tau_0}$ is only $C^{\infty}$--uniform on compact sets, the limit might decay asymptotically at a different rate than the solutions $\tu_{\tau_k}$. The proof of the theorem is organized in several steps gradually improving our situation. First, we want to get into the position where we can apply the implicit function theorem, either theorem \ref{main-implicit-function-theorem} or theorem \ref{extended-implicit-function-theorem}, to the limit solution $\tu_{\tau_0}$.\\

\underline{\bf First step: Compute the Maslov index of $\tu_{\tau_0}$:}\\
The solutions $\tu_{\tau}$, $\tau<\tau_0$ all satisfy $\lambda_{\pm}=\mp\frac{\pi}{2}$ and $\mu(\tu_{\tau})=0$. This means that we can find a smooth family $n_{\tau}$ of admissible normal vectors so that $n_{\tau}(s,1)\in T_{\tu_{\tau}(s,1)}(\{0\}\times{\mathcal D})$ for all $s\in {\bf R}$. In the paper \cite{part2} (lemma 3.2) we have constructed admissible normal vectors rather explicitly from the maps $\tilde{u}_{\tau}$. We will review this construction below. It is clear from this construction that the admissible normal vectors $n_{\tau_k}$ will converge in $C^{\infty}_{loc}$ to some $n_{\tau_0}$ which will not be normal to $\tu_{\tau_0}$ near the ends if the decay rate of $\tu_{\tau_0}$ is not $\pm\frac{\pi}{2}$. Otherwise the limit $n_{\tau_0}$ is also an admissible normal vector for $\tu_{\tau_0}$, and therefore $\mu(\tu_{\tau_0})=0$. Ends of $\tu_{\tau_0}$ with decay rate $\lambda_{\pm}=\mp\frac{\pi}{2}$ do not make any contribution to the Maslov index $\mu(\tu_{\tau_0})$. Therefore, we will assume that $\tu_{\tau_0}$ decays like $e^{-m\pi |s|}$, $m\in{\bf N}$, for large $|s|$ in at least one of the ends, and we will only discuss the end where the decay rate is not $\pm\frac{\pi}{2}$. Pick $\ve>0$ and $R>0$ so large that the asymptotic formula (theorem \ref{asymptotic-formula-theorem}) holds for $\tu_{\tau_0}(s,t)$ with $|s|\ge R$, and such that the remainder term in the asymptotic formula and its derivatives are no larger than $\ve$. We also want $R>0$ so large that $\pil\pas u_{\tau_0}(s,t)\neq 0$ for all $|s|\ge R$. It follows easily from the asymptotic formula (\ref{even-explicit-asymptotic-formula-for-us-v2}) applied to $\pas\tu_{\tau_0}$ that such an $R$ can be found. We will work in an open neighborhood $U$ of $\{0\}\times{\mathcal L}$ where the coordinates of proposition \ref{3.1.1.} can be used.
We use the following complex frame for the contact structure $\xi|_U$:
\[
e_1(\theta,x,y):=(0,1,0,-x)\ ,\ e_2(\theta,x,y):=-\tilde{J}(\theta,x,y)e_1=(0,0,1,0).
\]
Denote by $(\pi_{\lambda}\partial_su_{\tau_k})_1$ and $(\pi_{\lambda}\partial_su_{\tau_k})_2$ the components of $\pi_{\lambda}\partial_su_{\tau_k}$ along $e_1$ and $e_2$ respectively. Define now normal vectors to $\tu_{\tau_k}$ by
\begin{eqnarray}\label{def-of-nbar}
\bar{n}_k(s,t) & := &  \Big(-(\pi_{\lambda}\partial_su_{\tau_k})_1,\partial_sa_{\tau_k}\cdot e_1(u_{\tau_k})-(\pi_{\lambda}\partial_su_{\tau_k})_2X_{\lambda}(u_{\tau_k})+\nonumber\\
 & & +(\lambda(u_{\tau_k})\partial_su_{\tau_k})\cdot e_2(u_{\tau_k})\Big)(s,t),
\end{eqnarray}
and
\[
\bar{m}_k(s,t):=\tilde{J}(\tu_{\tau_k}(s,t))\bar{n}_k(s,t).
\]
Writing $\tu_{\tau_k}=(a,\theta,x,y)$ and 
\[
\Lambda=\lambda(u_{\tau_k})\pas u_{\tau_k}=\pas y+x\pas\theta
\]
we obtain
\begin{eqnarray*}
\pas\tu_{\tau_k} & = & \pas a\,\frac{\partial}{\partial\tau}+\Lambda\,X_{\lambda}+\pas\theta\,e_1+\pas x\,e_2\\
\pat\tu_{\tau_k} & = & -\Lambda\,\frac{\partial}{\partial\tau}+\pas a\,X_{\lambda}+\pas x\,e_1-\pas\theta\,e_2
\end{eqnarray*}
and
\begin{eqnarray*}
\bar{n}_k & = & -\pas\theta\,\frac{\partial}{\partial\tau}-\pas x\,X_{\lambda}+\pas a\,e_1+\Lambda\,e_2\\
\bar{m}_k & = & \pas x\,\frac{\partial}{\partial\tau}-\pas\theta\,X_{\lambda}+\Lambda\,e_1-\pas a\,e_2.
\end{eqnarray*}
Away from the boundary singular points the tangent space to $\{0\}\times{\mathcal D}$ at the point $\tu_{\tau_k}(s,1)$ is generated by the vectors
\[
e_1(u_{\tau_k})+x\,X_{\lambda}(u_{\tau_k})+q'(\theta)y\,e_2(u_{\tau_k})|_{(s,1)}\ \mbox{and}\ q(\theta)\,e_2(u_{\tau_k})+X_{\lambda}(u_{\tau_k})|_{(s,1)}.
\]
For $t=0$ we have
\[
\bar{n}_k(s,0)\,,\,\pas\tu_{\tau_k}(s,0)\in T_{\tu_{\tau_k}(s,0)}({\bf R}\times{\mathcal L}),
\]
and for $t=1$ the vector space $T_{\tu_{\tau_k}(s,1)}(\{0\}\times{\mathcal D})$ is generated by $\pas\tu_{\tau_k}(s,1)$ and
\begin{equation}\label{generator-of-boundary-condition}
q(\theta(s,1))\pat\tu_{\tau_k}(s,1)+\bar{m}_k(s,1)+y(s,1)(q'(\theta(s,1))-q^2(\theta(s,1)))\,\bar{n}_k(s,1).
\end{equation}
If $n_{\tau_k}$ is an admissible normal vector to $\tu_{\tau_k}$ then we can write it in the following form using the previously constructed normal vectors
\begin{equation}\label{compute-admissible-vector}
n_{\tau_k}(s,t)=\Big(\alpha_1\pas\tu_{\tau_k}+\alpha_2\pat\tu_{\tau_k}+\beta_1\bar{n}_k+\beta_2\bar{m}_k\Big)(s,t)
\end{equation}
with suitable smooth coefficient functions so that $\beta_1^2+\beta_2^2$ is never zero. By choosing $k$ large we may assume that for any partial derivative $D^{\alpha}$
\begin{equation}\label{approxi-at-R}
\sup_{0\le t\le 1}\left|D^{\alpha}(\tu_{\tau_k}-\tu_{\tau_0})(\pm R,t)\right|<\ve,
\end{equation}
in particular, for $s=\pm R$, we may use the asymptotic formula (\ref{even-explicit-asymptotic-formula-v2}) for $\tu_{\tau_0}$ also for describing $\tu_{\tau_k}$. Choose now $R'>R$ so large that the asymptotic formula (\ref{odd-explicit-asymptotic-formula-v2}) holds for $\tu_{\tau_k}(s,t)$, $|s|\ge R'$, and such that the remainder term and its derivatives are no larger than $\ve$. Because $n_{\tau_k}$ is admissible it must converge to $(0,\pm 1,0,0)$ as $s\rightarrow\pm\infty$. This imposes certain conditions on the coeficients in (\ref{compute-admissible-vector}). We write $(0,\theta^k_{\pm},0,0)=\lim_{s\rightarrow\pm\infty}\tu_{\tau_k}(s,t)$. Using formula (\ref{odd-explicit-asymptotic-formula-v2}) for $|s|\ge R'$ we get
\begin{eqnarray*}
\pas a(s,t) & = & \pm\kappa_{\pm}\rho(s)\cos\left(\frac{\pi t}{2}\right)+\rho(s)\ve(s,t)\\
\pas\theta(s,t) & = & \mp\kappa_{\pm}\rho(s)q(\theta^k_{\pm})\cos\left(\frac{\pi t}{2}\right)+\rho(s)\ve(s,t)\\
\pas x(s,t) & = & -\kappa_{\pm}\rho(s)q(\theta^k_{\pm})\sin\left(\frac{\pi t}{2}\right)+\rho(s)\ve(s,t)\\
\pas y(s,t) & = & -\kappa_{\pm}\rho(s)\sin\left(\frac{\pi t}{2}\right)+\rho(s)\ve(s,t)\\
\Lambda & = & -\kappa_{\pm}\rho(s)\sin\left(\frac{\pi t}{2}\right)+\rho(s)\ve(s,t)\\
q(\theta(s,t)) & = & q(\theta^k_{\pm})+\rho(s)\cdot\mbox{"something bounded"'},
\end{eqnarray*}
where $\kappa_{\pm}>0$ are some constants and $\rho(s)=e^{\int_{s_0}^s\alpha_{\pm}(\tau)d\tau}$ is the function obtained by applying theorem \ref{asymptotic-formula-theorem} to $\tu_{\tau_k}$. As usual, we denote by $\ve(s,t)$ a smooth function which converges to zero with all its derivatives uniformly in $t$ as $|s|\rightarrow\infty$. Using the above asymptotic formulae, the condition $n_{\tau_k}(s,t)\rightarrow (0,\pm 1,0,0)$ as $s\rightarrow\pm\infty$ has the following implication on the coefficients in (\ref{compute-admissible-vector}):
\begin{eqnarray}\label{formula-for-admissible-vector-2}
n_{\tau_k}(s,t) & = & \frac{1}{\kappa_{\pm}\rho(s)(1+q^2(\theta^k_{\pm}))}\Big[-q(\theta^k_{\pm})\cos\left(\frac{\pi t}{2}\right)\pas \tu_{\tau_k}(s,t)\mp\nonumber\\
 & & \mp q(\theta^k_{\pm})\sin\left(\frac{\pi t}{2}\right)\pat\tu_{\tau_k}(s,t)+\cos\left(\frac{\pi t}{2}\right)\bar{n}_k(s,t)\mp\\
  & & \mp\sin\left(\frac{\pi t}{2}\right)\bar{m}_k(s,t)+\rho(s)\ve(s,t)\Big]\nonumber.
\end{eqnarray}
Without the smaller order term, equation (\ref{formula-for-admissible-vector-2}) would not be correct because the boundary condition for $t=1$ would not be satisfied. We may write equation (\ref{formula-for-admissible-vector-2}) as follows
\begin{eqnarray*}
n_{\tau_k}(s,t) & = & \frac{1}{\kappa_{\pm}\rho(s)(1+q^2(\theta^k_{\pm}))}\Big[-q(\theta^k_{\pm})\cos\left(\frac{\pi t}{2}\right)\pas \tu_{\tau_k}(s,t)\mp\\
 & & \mp q(\theta(s,t))\sin\left(\frac{\pi t}{2}\right)\pat\tu_{\tau_k}(s,t)+\cos\left(\frac{\pi t}{2}\right)\bar{n}_k(s,t)\mp\\
  & & \mp\sin\left(\frac{\pi t}{2}\right)\bar{m}_k(s,t)+\rho(s)\ve(s,t)\mp\\
  & & \mp\sin\left(\frac{\pi t}{2}\right) y(s,t)(q'(\theta(s,t))-q^2(\theta(s,t)))\bar{n}_k(s,t)\Big].
\end{eqnarray*}
because the expressions
\[
\mp q(\theta^k_{\pm})\pat\tu_{\tau_k}(s,t)\mp\bar{m}_k(s,t)
\]
and
\[
\mp q(\theta(s,t))\pat\tu_{\tau_k}(s,t)\mp\bar{m}_k(s,t)\mp y(s,t)(q'(\theta(s,t))-q^2(\theta(s,t)))\bar{n}_k(s,t)
\]
differ only by a term of the form $\rho(s)\ve(s,t)$. We may now change the admissible normal vector $n_{\tau_k}$ by removing the remainder term $\rho(s)\ve(s,t)$ with a smooth cut--off function for $|s|\ge R$. In the same way we may we may also replace all the terms $q(\theta^k_{\pm})$ with $q(\theta(s,t))$. We will keep the notation $n_{\tau_k}$ for simplicity, and we still have an admissible normal vector to $\tu_{\tau_k}$ which satisfies both boundary conditions for $t=0$ and $t=1$. As we remarked earlier, the vectors $n_{\tau_k}$ converge in $C^{\infty}([-R,R]\times [0,1])$ to some normal vector $n_{\tau_0}$ to $\tu_{\tau_0}|_{[-R,R]\times[0,1]}$ which also satisfies the boundary conditions $n_{\tau_0}(s,0)\in T_{\tu_{\tau_0}(s,0)}({\bf R}\times{\mathcal L})$ and $n_{\tau_0}(s,1)\in T_{\tu_{\tau_0}(s,1)}(\{0\}\times{\mathcal D})$. In local coordinates $n_{\tau_0}$ is given by the formula
\begin{eqnarray}\label{limit-normal-vector}
n_{\tau_0}(s,t) & = & \frac{1}{\kappa_{\pm}\rho_0(s)(1+q^2(\theta_0(s,t)))}\Big[-q(\theta_0(s,t))\cos\left(\frac{\pi t}{2}\right)\pas \tu_{\tau_0}(s,t)\mp\nonumber\\
 & & \mp q(\theta_0(s,t))\sin\left(\frac{\pi t}{2}\right)\pat\tu_{\tau_0}(s,t)+\cos\left(\frac{\pi t}{2}\right)\bar{n}_0(s,t)\mp\nonumber\\
  & & \mp\sin\left(\frac{\pi t}{2}\right)\bar{m}_0(s,t)\mp\\
  & & \mp\sin\left(\frac{\pi t}{2}\right) y_0(s,t)(q'(\theta_0(s,t))-q^2(\theta_0(s,t)))\bar{n}_0(s,t)\Big],\nonumber
\end{eqnarray}
where the subscript '0' indicates that we use $\tu_{\tau_0}$ for evaluating the formula instead of $\tu_{\tau_k}$. We can now extend the normal vector $n_{\tau_0}$ to the whole strip $S$ simply by using (\ref{limit-normal-vector}) above for all $s\in{\bf R}$. It will not be an admissible normal vector, the behavior at infinity is not the same as in definition \ref{admissible-normal-vector}, otherwise we would have $\mu(\tu_{\tau_0})=0$. Using the asymptotic formula (\ref{even-explicit-asymptotic-formula-v2}) for $\tu_{\tau_0}$ we will determine its winding behavior at infinity which enables us to compute $\mu(\tu_{\tau_0})$. Since we are only interested in the limits
\[
\tilde{n}_{\pm}(t)=\lim_{s\rightarrow\pm\infty}n_{\tau_0}(s,t)
\]
we can ignore the remainder terms in the asymptotic formulae for $\tu_{\tau_0}$. We may also neglect the term
\[
\frac{1}{\rho_0(s)}(\sin\left(\frac{\pi t}{2}\right) y_0(s,t)(q'(\theta_0(s,t))-q^2(\theta_0(s,t)))\bar{n}_0(s,t))
\]
since it tends to zero in the limit $s\rightarrow\pm\infty$. Using (\ref{limit-normal-vector}) and
\begin{eqnarray*}
\pas\tu_{\tau_0}(s,t) & = & \kappa\rho_0(s)\,(0,\mp\cos(m_{\pm}\pi t), -\sin(m_{\pm}\pi t),0)\\
\pat\tu_{\tau_0}(s,t) & = & \kappa\rho_0(s)\,(0,-\sin(m_{\pm}\pi t),\pm\cos(m_{\pm}\pi t),0)\\
\bar{n}_0(s,t) & = & \kappa\rho_0(s)\,(\pm\cos(m_{\pm}\pi t),0,0,\sin(m_{\pm}\pi t))\\
\bar{m}_0(s,t) & = & \kappa\rho_0(s)\,(-\sin(m_{\pm}\pi t),0,0,\pm\cos(m_{\pm}\pi t))
\end{eqnarray*}
with some nonzero constant $\kappa$ we obtain
\begin{eqnarray*}
\tilde{n}_{\pm}(t) & = & c\,\left(\pm\cos\left(\frac{2m_{\pm}-1}{2}\pi t\right),\pm q(\theta^0_{\pm})\cos\left(\frac{2m_{\pm}-1}{2}\pi t\right),\right.\\
 & & \left.q(\theta^0_{\pm})\sin\left(\frac{2m_{\pm}-1}{2}\pi t\right),\sin\left(\frac{2m_{\pm}-1}{2}\pi t\right)\right)
\end{eqnarray*}
(with some nonzero constant $c$). In definition \ref{admissible-normal-vector} we have used the coordinates (\ref{make-boundary-conditions-flat}). In the coordinates given by proposition \ref{3.1.1.} the corresponding vector $n_{\pm}(t)$ is given by
\[
(\cos\frac{\pi t}{2},-q(\theta^0_{\pm})\cos\frac{\pi t}{2},\mp q(\theta^0_{\pm})\sin\frac{\pi t}{2},\mp\sin\frac{\pi t}{2})
\]
and
\[
m_{\pm}(t)=\tilde{J}(0,\theta^0_{\pm},0,0)n_{\pm}(t)=(\pm \sin\frac{\pi t}{2}, \mp q(\theta^0_{\pm})\sin\frac{\pi t}{2}, q(\theta^0_{\pm})\cos\frac{\pi t}{2},\cos\frac{\pi t}{2}).
\]
With 
\begin{eqnarray*}
e(t) & := & (0,\mp\cos(m_{\pm}\pi t), -\sin(m_{\pm}\pi t),0)\\
f(t) & := & (0,-\sin(m_{\pm}\pi t),\pm\cos(m_{\pm}\pi t),0)
\end{eqnarray*}
we will now compute coefficients $\alpha^{\pm}_1,\alpha^{\pm}_2,\beta^{\pm}_1,\beta^{\pm}_2$ depending on $t$ such that
\[
\tilde{n}_{\pm}(t)=\beta^{\pm}_1(t)n_{\pm}(t)+\beta^{\pm}_2(t)m_{\pm}(t)+\alpha^{\pm}_1(t)e(t)+\alpha^{\pm}_2(t)f(t).
\]
The relationship with the Maslov index $\mu(\tu_{\tau_0})$ is now the following: Writing 
\[
\frac{\beta^{\pm}_1+i\beta^{\pm}_2}{|\beta^{\pm}_1+i\beta^{\pm}_2|}(t)=e^{i\psi_{\pm}(t)}
\]
for suitable smooth functions $\psi_{\pm}$ we have
\[
\mu(\tu_{\tau_0})=\frac{1}{\pi}(\psi_+(0)-\psi_+(1)+\psi_-(1)-\psi_-(0)).
\]
The computation of the coefficients yields
\[
\beta^{\pm}_1(t)+i\beta^{\pm}_2(t)=\pm c\,e^{\pm\,im_{\pm}\pi t}
\]
so that 
\[
\mu(\tu_{\tau_0})\ =\ -m_--m_+,
\]
i.e. for each end that has a decay rate of $m\pi$, $m\in{\bf Z}$ the Maslov--index changes by $-|m|$. The crucial remark is that the solution $\tu_{\tau_0}$ satisfies now one of the assumptions of the implicit function theorem (theorem \ref{extended-implicit-function-theorem}). In the next step we will show that $\tu_{\tau_0}$ is embedded.\\

\underline{\bf Second step: Show that $\tu_{\tau_0}$ is an embedding:}\\
In order to apply theorem \ref{extended-implicit-function-theorem} to the solution $\tu_{\tau_0}$ we have to make sure that it is immersed. It follows from the asymptotic formula, theorem \ref{asymptotic-formula-theorem}, and its versions in local coordinates that there is $R>0$ such that $\pas\tu_{\tau_0}(s,t)\neq 0$ whenever $|s|\ge R$, regardless of the decay behavior of $\tu_{\tau_0}$. On the other hand, it follows from proposition \ref{maximum-principle} that $\pas\tu_{\tau_0}(s,1)\neq 0$ for all $s\in{\bf R}$ as well. The idea is now to use a result about positivity of intersections of pseudoholomorphic curves in four dimensional almost complex manifolds: Assume that $u$ is a non-constant pseudoholomorphic disk in an almost complex manifold $(W,J)$ where $\hbox{dim}\,W=4$:
\begin{eqnarray*}
&u:D\rightarrow W&\\
&\pas u +J(u)\pat u\,=\,0&
\end{eqnarray*}
We say that $u$ is an embedding near the boundary if the following holds: There exists a small annulus around the boundary $A_{\varepsilon}$,
$$
A_{\varepsilon} =\{z\in D |\ 1-\varepsilon \leq |z|\leq 1\}
$$
such that 
\begin{eqnarray*}\label{embeddingboundary}
&u|A_{\varepsilon}\ \hbox{is an embedding}&\\
&u^{-1}(u(A_{\varepsilon})) =A_{\varepsilon}.&
\end{eqnarray*}
For such an embedding at the boundary one can define a self-intersection index $I(u)\in{\bf Z}$ (see \cite{AH}, \cite{Gromov},\cite{McDuff2},\cite{Sikorav}) which has the following properties:
\begin{itemize}
\item If $u_{\tau}$ is a smooth family of pseudoholomorphic disks which are embeddings at the boundary then the intersection indices $I(u_{\tau})$ are independent of $\tau$,
\item $I(u)=0$ if and only if $u$ has no singularities and no self--intersections.
\end{itemize}
It follows from proposition \ref{reflection-2} that we may treat boundary points on ${\bf R}\times\{0\}$ like interior points since the solutions $\tu_{\tau_0},\tu_{\tau_k}$ can be locally reflected at the boundary. We recall that $\tu_{\tau_0}$ is approximated in $C^{\infty}([-R,R]\times[0,1])$ by the sequence $\tu_{\tau_k}$ which are all embedded solutions. It is well known (see \cite{AH},\cite{McDuff2}) that the following alternative holds for each point $z\in S$: Either there is some $\delta>0$ such that $\tu_{\tau_0}|_{B_{\delta}(z)\backslash\{z\}}$ is an embedding or there is a biholomorphic map $\phi:U\rightarrow V$ between neighborhoods of $z$ such that $\phi(z')\neq z'$ for some $z'\in U$ and $\tu_{\tau_0}|_V\circ \phi=\tu_{\tau_0}|_U$. We claim that there is a point $z_0\in{\bf R}\times\{1\}$ such that 
\[
\tu_{\tau_0}^{-1}(\tu_{\tau_0}(z_0))=\{z_0\}.
\]
Let us first show how to finish the proof of the immersion property assuming that the claim is correct. Consider the set
\[
{\mathcal S}:=\{z\in S\,|\,\pas\tu_{\tau_0}(z)\neq 0\},
\]
and the set ${\mathcal S}'$ consisting of all $z\in {\mathcal S}$ such that there is $z'\neq z$ and sequences $z_k\rightarrow z$, $z'_k\rightarrow z'$ (of course we also assume $z_k\neq z$, $z'_k\neq z'$) with $\tu_{\tau_0}(z'_k)=\tu_{\tau_0}(z_k)$. The set ${\mathcal S}'$ is closed in ${\mathcal S}$, but it is also open by the Similarity Principle (see appendix \ref{similarity-principle}). In fact, the Similarity Principle implies that intersection points between pseudoholomorphic curves can only accumulate in a point which is critical on both curves unless the images of the two curves coincide. Therefore either ${\mathcal S}'=\emptyset$ or ${\mathcal S}'={\mathcal S}$. The latter alternative cannot hold because of the claim above and because $\pas\tu_{\tau_0}(s,0)\neq 0$ for all $s\in{\bf R}$ by proposition \ref{maximum-principle}. Recalling that $\tu_{\tau_0}$ is not constant we conclude from the Similarity Principle that the critical points of $\tu_{\tau_0}$ are isolated. If for given $w\in S$ we could find a nontrivial biholomorphic map $\phi$ between neighborhoods of $w$ such that $\tu_{\tau_0}\circ\phi=\tu_{\tau_0}$ then ${\mathcal S}'$ would not be empty: Just pick any point $z$ near $w$ which is not critical, a sequence $z_k$ converging to it, and take $z'=\phi(z)$, $z'_k=\phi(z_k)$.
Hence for any point $z\in S$ we can find $\delta>0$ such that $\tu_{\tau_0}$ restricted to the punctured neighborhood $B_{\delta}(z)\backslash\{z\}$ is an embedding. This means that the self--intersection index of $\tu_{\tau_0}|_{B_{\delta}(z)}$ is well--defined. The maps $\tu_{\tau_k}|_{B_{\delta}(z)}$ are all embeddings. Hence they have zero self--intersection index and so does $\tu_{\tau_0}|_{B_{\delta}(z)}$ for any $z\in S$. Therefore, $\tu_{\tau_0}$ must be an immersion provided the claim we made earlier is correct.\\
We still have to show that there is a point $z_0\in{\bf R}\times\{1\}$ such that 
\[
\tu_{\tau_0}^{-1}(\tu_{\tau_0}(z_0))=\{z_0\}.
\]
We are going to show more, we will actually prove that the curve
\[
s\longmapsto u_{\tau_0}(s,1)\in{\mathcal D}^{\ast}
\]
has no self--intersections at all. We argue by contradiction, and we assume that there are $z_0,z_1\in{\bf R}\times\{1\}$ with $z_0\neq z_1$ and $\tu_{\tau_0}(z_0)=\tu_{\tau_0}(z_1)$. By proposition \ref{maximum-principle} and corollary \ref{local-reflection} we may reflect $\tu_{\tau_0}$ near the boundary points $z_0$ and $z_1$. Locally near $\tu_{\tau_0}(z_0)$ we are in the following situation: We have two pseudoholomorphic disks $u,v:D\rightarrow{\bf C}^2$ with respect to some almost complex structure $\bar{J}$ on ${\bf C}^2$ with $\bar{J}(0)=i$. In addition, we have $u(0)=v(0)$, and $0\in D$ is not a critical point for any of the maps $u$ and $v$. This implies that $0$ is an isolated intersection point, hence we may assume that $u-v$ is not zero on $D\backslash\{0\}$. The disks $u$ and $v$ are approximated by disks $u_k$, $v_k$ with $\inf_D|u_k-v_k|>0$. We obtain a contradiction since the algebraic intersection number of $u_k(D)$ and $v_k(D)$ is zero for all $k$ while it is at least one for $u(D)$ and $v(D)$ (see \cite{AH}, \cite{McDuff2}). This completes the proof of the claim. We point out that this argument actually shows that $\tu_{\tau_0}$ can only have a self--intersection in points $z_0,z_1$ which are both critical, i.e. where $\pas\tu_{\tau_0}(z_0)=\pas\tu_{\tau_0}(z_1)=0$.\\
Let us summarize and proceed to the embedding property: We know by proposition \ref{maximum-principle} that $\pas\tu_{\tau_0}(s,1)\neq 0$ for all $s\in{\bf R}$. The argument outlined above then shows that the curve $s\mapsto u_{\tau_0}(s,1)$ does not have any self--intersections which in turn implies the claim. On the other hand, we then know that $\tu_{\tau_0}$ is immersed. Since self--intersection points of $\tu_{\tau_0}$ can only occur in critical points, we know that there are none. Because there is also no point $z\in S$ with $\tu_{\tau_0}(z)=\lim_{s\rightarrow\pm\infty}\tu_{\tau_0}(s,t)$ we conclude that $\tu_{\tau_0}$ is an embedding. Indeed, if we had $\tu_{\tau_0}(z)=\lim_{s\rightarrow\pm\infty}\tu_{\tau_0}(s,t)\in\{0\}\times{\mathcal L}$ then proposition \ref{maximum-principle} would imply $z\in{\bf R}\times \{1\}$. Because the curve $s\mapsto u_{\tau_0}(s,1)\in{\mathcal D}$ is always transverse to the characteristic foliation, it can never hit the boundary ${\mathcal L}$, i.e. there is no such point $z$. 

{\bf \underline{Third step:} Show that $\tu_{\tau_0}$ has the same rate of decay as $\tu_{\tau_k}$, that its Maslov index vanishes and that it is the limit of the maps $\tu_{\tau}$ as $\tau\nearrow\tau_0$:}\\

Until now, the convergence $\tu_{\tau_k}\rightarrow\tu_{\tau_0}$ is uniform only on compact sets. Nevertheless, we have succeeded to verify all the assumptions of the implicit function theorem (theorem \ref{extended-implicit-function-theorem}). If the decay rates of $\tu_{\tau_0}$ are already $\lambda_{\pm}=\mp\frac{\pi}{2}$ then we also have $\mu(\tu_{\tau_0})=0$ and we can apply the original implicit function theorem from the paper \cite{part2} (theorem \ref{main-implicit-function-theorem}). In any case, there is an $N$--dimensional family of solutions $(\tilde{v}_{\sigma})_{\sigma\in{\bf R}^N}$ with $\tilde{v}_0=\tu_{\tau_0}$. The implicit function theorems also imply that the family $\tilde{v}_{\sigma}$ is the only family of solutions close to $\tu_{\tau_0}$ in the following sense: Any other solution whose image lies in the neighborhood $U$ of $\tu_{\tau_0}(S)$ used in the proof of theorems \ref{extended-implicit-function-theorem} or \ref{main-implicit-function-theorem}, must be one of the solutions $\tilde{v}_{\sigma}$ up to parametrization. It is clear that this is not sufficient for our purpose because at this point it may be possible that $\tu_{\tau_k}\rightarrow \tu_{\tau_0}$ in $C^{\infty}_{loc}$, but the image of $\tu_{\tau_k}$ is not contained in a neighborhood of $\tu_{\tau_0}(S)$. This phenomenon occurs in Morse theory and Floer homology if trajectories converge to a broken trajectory. Our situation is different because of our intersection result, theorem \ref{interior-and-boundary-intersections}. We note that the families $\tilde{v}_{\sigma}$ and $\tu_{\tau}$ do not intersect for small $\tau$ and $|\sigma|$. On the other hand they have to intersect later, say for some $\sigma'$ and $\tau'$. The union $V$ of the curves $(\tilde{v}_{\sigma}((-R,R)\times\{1\}))_{\sigma\in{\bf R}^N}$ is a neighborhood of the curve $\tu_{\tau_0}((-R,R)\times\{1\})$ on the Seifert surface ${\mathcal D}$ , and for large $k$ the curves $\tu_{\tau_k}((-R,R)\times\{1\})$ have to enter $V$. Of course, we took advantage of the two--dimensional situation here. Since there is no isolated first intersection between the family $\tilde{v}_{\sigma}$ and the family $\tu_{\tau}$, the images of $\tilde{v}_{\sigma'}$ and $\tu_{\tau'}$ agree. In this case $\tilde{v}_{\sigma'}$ can not agree with the image of some sort of Schwarz reflection of $\tu_{\tau'}$ because both are close to $\tau_{\tau_0}$ on compact sets. This implies that the image of the solution $\tu_{\tau'}$ is in fact close to the image of $\tu_{\tau_0}$, and the families $\tilde{v}_{\sigma}$ and $\tu_{\tau}$ actually are just one family. Hence the decay rates of $\tu_{\tau_0}$ and all the $\tilde{v}_{\sigma}$ are $|\lambda_{\pm}|=\frac{\pi}{2}$, the Maslov indices are all zero, the family is in fact one--dimensional ($N=1$), and it is produced by the original implicit function theorem, theorem \ref{main-implicit-function-theorem}. Therefore, we obtain the same limit $\tu_{\tau_0}$ for all sequences $\tau_k\nearrow\tau_0$ and the convergence is in $C^{\infty}(S)$, not just in $C^{\infty}_{loc}$. This completes the proof of theorem \ref{compactness-result}. \qed 

\begin{appendix}
\section{The Similarity principle}\label{similarity-principle}

In this appendix, we review the similarity principle in its original version and also a version near boundary points. 
For $ 2 < p < \infty$ we denote by $V^p$ the Banach space
consisting of all $u \in W^{1,p} (D, {\bf C}^n)$ satisfying $u
(\partial D) \subset {\bf R}^n$, where $D\subset{\bf C}$ is the
unit disk and let $$ \bar{\partial}
: V^p \rightarrow L^p $$ be the standard Cauchy Riemann operator $$u
\longmapsto \frac{\partial u}{\partial s} + i \frac{\partial
u}{\partial t}. $$The operator $\bar{\partial} : V^p \rightarrow L^p(D,{\bf C}^n)$
is surjective and Fredholm of index $n$ with kernel being
the constants in ${\bf R}^n$.

\begin{theorem} \label{A4.theorem 1}
Assume $A \in L^{\infty} (D, {\cal L}_{\bf R} ({\bf C}^n))$ ,
$2<p<\infty$ and $ w \in W^{1,p}_{loc}(\stackrel{\circ}{D},{\bf C}^n)$. Let $w$ be a
solution of
\begin{eqnarray*}
  \bar{\partial} w + Aw & = & 0 \qquad \mbox{in} \qquad \stackrel{\circ}{D}\\
  w (0)                 & = & 0.
\end{eqnarray*}
Then there exists $$\Phi \in \bigcap_{2 < q < \infty} W^{1,q}
(D, {\cal L}_{\bf C} ({\bf C}^n))$$ with $$\Phi (0)
= \Id\ \ , \ \ \Phi (z) \in \GL ({\bf C}^n)$$ and a map $f : D\rightarrow {\bf C}^n$ with $f(0)=0$ such
that for $z\in D$ 
$$
  w (z) = \Phi (z) f (z).
$$  
Moreover, if $0<\varepsilon\le 1$ is a sufficiently small number and
$D_{\varepsilon}\subset{\bf C}$ the disk of radius $\varepsilon$ then $f$ is holomorphic on $D_{\varepsilon}$.
\end{theorem}

{\bf Proof of theorem \ref{A4.theorem 1}:} See \cite{AH} or \cite{HZbook}\qed 

Next we consider a boundary version of the similarity principle. Let $$D^+ := \{z \in D \mid \mbox{Im}
(z) \geq 0\}$$

\begin{theorem} \label{A4.theorem 2}\label{similarity-principle-bdry}
Assume $A \in L^{\infty} (D^+ , {\cal L}_{\bf R} ({\bf C}^n))$ and $w \in
W_{loc}^{1,p} (D^+ , {\bf C}^n)$, $2 < p < \infty$, satisfying
\begin{eqnarray*}
  \bar{\partial} w + A w & = & 0 \qquad \mbox{on} \qquad \stackrel{\circ}{D}^+\\
  w ((-1,1)) & \subset & {\bf R}^n, \qquad w(0) = 0.
\end{eqnarray*}
Then there exists $\Phi \in \bigcap_{2 < q < \infty} W^{1,q}
(D^+, {\cal L}_{\bf C} ({\bf C}^n))$ with $$
     \Phi(z)\in\GL({\bf C}^n)\ \ ,\ \ \Phi(0)=\mbox{Id}
$$
$$
     \Phi (z) \in \GL({\bf R}^n)\subset{\cal L} ({\bf R}^n) \qquad \mbox{for} \quad z \in
     (-1, 1)
$$
and a map $f : D^+ \rightarrow {\bf C}^n$ with
$$
  f (z) \in {\bf R}^n \qquad \mbox{for} \quad z \in (-1,1)\ ,\ f(0)=0,
$$
holomorphic on some smaller half--disk $D^+_{\varepsilon}$, such that
$$
  w (z) = \Phi (z) f (z).
$$
\end{theorem}

{\bf Proof:}\\ This result can be reduced to Theorem
\ref{A4.theorem 1}. Extend $A$ to a map in $L^{\infty} (D, {\cal
L}_{\bf R} ({\bf C}^n))$ by $$
  A (z) = \overline{A (\bar{z})}\ \mbox{if  Im}(z)<0,
$$
where `` $\overline{\phantom{zzzz}} $ '' means replacing all coefficients by
the complex conjugate ones. Extend $w$ similarly by
$$
  w(z) = \overline{w(\bar{z})}.
$$
Then $w \in W^{1,p} (D, {\bf C}^n)$ as one verifies easily. Now apply
 theorem \ref{A4.theorem 1} and find
$$
  w (z) = \Phi (z) f (z),
$$ where $z$ lies in some disk $D_{\varepsilon}$, and it turns out that  $$
  \Phi ((-1,1)) \subset {\cal L} ({\bf R}^n)
$$
and consequently
$$
  f ((-\varepsilon, \varepsilon)) \subset {\bf R}^n.
$$ \qed

{\bf Remark:}\\
There is also a parameterized version of the Similarity Principles: If $A_{\tau}$ is a continuous path in $L^{\infty}(D,{\mathcal L}_{{\bf R}}({\bf C}^n))$ with $A_0=0$ and if $w_{\tau}$ is a continuous family of solutions of $\bar{\partial}w_{\tau}+A_{\tau}w_{\tau}=0$ then $w_{\tau}=\Phi_{\tau}\sigma_{\tau}$ with $\Phi_{\tau},\sigma_{\tau}\in C^0(D)$ depending continuously on $\tau$ as well.  The maps $\Phi_{\tau}$ converge in $C^0(D)$ to the identity matrix. The important fact is that the path of operators $\Phi\mapsto (\bar{\partial}\Phi+A_{\tau}\Phi,\Phi(1))$ and the corresponding path of the inverses are continuous in $\tau$ with respect to the operator norm. 

\end{appendix}

\end{document}